\documentclass[12pt]{amsart}

\usepackage{amsfonts,amsmath,amssymb,amsthm,mathrsfs,url,xcolor,latexsym,rawfonts,graphicx}
\usepackage{hyperref}
\usepackage{enumerate} 
\usepackage[]{caption2}
\usepackage{accents}
\usepackage{amssymb,amsmath}
\usepackage{amsbsy}
\usepackage{cite}
\usepackage{float}
\usepackage{graphicx}

\theoremstyle{plain}

   \newcommand{\beqn}{\begin{eqnarray}}
   \newcommand{\eeqn}{\end{eqnarray}}
   \newcommand{\beqs}{\begin{eqnarray*}}
   \newcommand{\eeqs}{\end{eqnarray*}}
   \newcommand{\ban}{\begin{eqnarray*}}
   \newcommand{\nan}{\end{eqnarray*}}
   \newcommand{\beq}{\begin{equation}}
   \newcommand{\eeq}{\end{equation}}

\renewcommand{\det}{\mbox{det}}

\numberwithin{equation}{section}

  \numberwithin{equation}{section}
  \numberwithin{figure}{section}

\begin{document}

\title[On the $L_\MakeLowercase{p}$ Gaussian Minkowski problem]
{\textbf{On the $L_\MakeLowercase{p}$ Gaussian Minkowski problem}}

\author[Yibin Feng]
{Yibin Feng}
\address{Yibin Feng, School of Mathematical Sciences,
University of Science and Technology of China,
Hefei, 230026, P.R. China.}
\email{fybt1894@ustc.edu.cn}

\author[Shengnan Hu]
{Shengnan Hu}
\address
	{Shengnan Hu, School of Mathematical Sciences,
University of Science and Technology of China,
Hefei, 230026, P.R. China.}
\email{helenhsn@mail.ustc.edu.cn}

\author[Lei Xu]
{Lei Xu}
\address
{Lei Xu, School of Mathematical Sciences,
University of Science and Technology of China,
Hefei, 230026, P.R. China.}
\email{xlsx@mail.ustc.edu.cn}

\subjclass[2010]{35J60, 52A40, 52A38.}

\keywords{Gaussian volume, $L_p$ Gaussian Minkowski problem, Monge-Amp\`{e}re equation, $L_p$ Gaussian surface area measure, convex body.}

\thanks{This work was supported by the China Postdoctoral Science Foundation (Grant No.2020TQ0315). 
 }


\date{\today}

\dedicatory{}

\begin{abstract}  
Existence of symmetric (resp. asymmetric) solutions to the $L_p$ Gaussian Minkowski problem for $p\leq 0$ (resp. $p\geq 1$) will be provided. Moreover, existence and uniqueness of smooth solutions to the problem for $p>n$ will also be proved without the restriction that the Gaussian volumes of convex bodies are not less than one-second.
\end{abstract}

\maketitle

\baselineskip=16.4pt
\parskip=3pt

\section{\bf Introduction}
One of the most important problems in convex geometric analysis is the Minkowski problem. The classical Minkowski problem asks for necessary and sufficient conditions on a given measure so that it is the surface area measure of a convex body. The existence of its solutions is due to Minkowski in the case of discrete measures and to Aleksandrov in the case of general measures. The study of the problem strong influence on fully nonlinear PDEs; see for example Lewy\cite{1}, Nirenberg\cite{2}, Cheng and Yau \cite{3}, Pogorelov \cite{4} and Caffarelli\cite{5}.

In \cite{6}, Lutwak showed that there is an $L_p$ analogue of the surface area measure and posed the associated $L_p$ Minkowski problem which has the classical Minkowski problem when $p=1$. The $L_p$ Minkowski problem for $p>1$ was solved by Lutwak \cite{6} for symmetric measures, and by Chou and Wang \cite{7} and Hug, Lutwak, Yang and Zhang \cite{8} for general measures. The problem for $p<1$ is much more complicated and contains two unsolved cases of the logarithmic Minkowski problem ($p=0$), see e.g., \cite{9, 10, 11, 16, 17, 19, 20, 21, 23, 32} and the centro-affine Minkowski problem ($p=-n$), see e.g., \cite{7, 24, 26, 27, 28, 29}. Important contributions to solving various cases of the $L_p$ Minkowski problem have a number of works; see say \cite{30, 31, 36, 37, 39, 40, 42, 49}.

Many different types of variations of the Minkowski problem have been introduced and extensively studied; see for example\cite{50, 51, 52, 53, 54, 55, 56, 57, 58, 59 ,60, 61, 62, 63, 64, 65, 66, 67, 68, 69, 70, 1q, 2q}. There is a particularly interesting Minkowski type problem, the Gaussian Minkowski problem, which was first introduced by Huang, Xi and Zhao \cite{71}. The origin of this new problem is the replacement of the surface area measure by the Gaussian surface area measure. The problem concerns both the existence and the uniqueness questions. The existence problem is to find necessary and sufficient conditions so that a given measure on the unit sphere is the Gaussian surface area measure of a convex body. The uniqueness question asks to what extent a convex body is uniquely determined by its Gaussian surface area measure.

Because of the lack of translation invariance and homogeneity for the Gaussian surface area measure, the Gaussian Minkowski problem is quite different from the classical Minkowski problem. Establishing the existence of solutions to the Gaussian Minkowski problem is more difficult and complicated than that to the Minkowski problem. Huang, Xi and Zhao \cite{71} considered the even case for the existence. That is, when the given measure $\mu$ is an even measure on the unit sphere that is not concentrated on any closed hemisphere and $|\mu|<\frac{1}{\sqrt{2\pi}}$, there exists an origin-symmetric convex body $K$ with larger than one-second Gaussian volume such that the Gaussian surface area measure of $K$ is just $\mu$. It is well-known that the density of the Gaussian surface area measure of a centered ball of radius $r$ is $\frac{1}{(\sqrt{2\pi})^n}e^{-\frac{r^2}{2}}r^{n-1}$. When $r$ approaches $0$ and $\infty$, $e^{-\frac{r^2}{2}}r^{n-1}\rightarrow 0$. This means that solutions to the Gaussian Minkowski problem are not unique without additional conditions. In \cite{71}, the uniqueness of solutions to the problem was established when restricted to convex bodies with the Gaussian volume not less than $\frac{1}{2}$. Since the seminal work of Huang, Xi and Zhao \cite{71}, the Gaussian Minkowski problem quickly became the center of attention; see for example \cite{72, 73, 74}. In particular, the existence results of non-symmetric solutions were very recently obtained in \cite{3q}.

Maybe as important as the $L_p$ Minkowski problem is the $L_p$ Gaussian Minkowski problem, posed in \cite{73}, prescribing the $L_p$ Gaussian surface area measure. The problem for $p=1$ is the classical Gaussian Minkowski problem
\cite{71}. In \cite{73}, the $L_p$ Gaussian Minkowski problem was solved in the general (not necessarily even) case when $p>0$, and in the even case when $p\geq 1$. The uniqueness of solutions to the problem  for $p\geq 1$ was also provd in  \cite{73}.

In this paper, we study the existence of symmetric (resp. asymmetric) solutions to the $L_p$ Gaussian Minkowski problem for $p\leq 0$ (resp. $p\geq 1$). In addition, the existence and uniqueness of solutions for $p>n$ in the smooth category will be obtained without the restriction that convex bodies have the Gaussian volume is great or equal to one-second. The related results are in detail stated as follows.

Let $\mathbb{R}^n$ be the $n$-dimensional Euclidean space. The unit sphere in $\mathbb{R}^n$ is denoted by $\mathbb{S}^{n-1}$. A convex body in $\mathbb{R}^n$ is a compact convex set with non-empty interior. Denote by $\mathcal{K}_o^n$ the set of all convex bodies that contain the origin in their interiors in $\mathbb{R}^n$ and by $\mathcal{K}_e^n$ the set of all origin-symmetric convex bodies in $\mathbb{R}^n$. The standard inner product of the vectors $x, y\in \mathbb{R}^n$ is denoted by $\langle x, y\rangle$. We write $|x|=\sqrt{\langle x, x\rangle}$ and use  $\mathcal{H}^{k}$ to denote the $k$-dimensional Hausdorff measure.

 The Gaussian volume, denoted by $\gamma_n$, of a convex body $K\in\mathcal{K}_o^n$ is defined by
 \begin{eqnarray*}
 	\gamma_n(K)=\frac{1}{\left(\sqrt{2\pi}\right)^n}\int_Ke^{-\frac{|x|^2}{2}}dx.
 \end{eqnarray*}
The $L_p$ Gaussian surface area measure $S_{p, \gamma_n}(K, \cdot)$, on $\mathbb{S}^{n-1}$, of  $K \in\mathcal{K}_o^n$ is uniquely determined by the following variational formula, discovered in \cite{73}, 
 \begin{eqnarray*}
 	\lim_{t\rightarrow 0}\dfrac{\gamma_n([h_t])-\gamma_n(K)}{t}=\frac{1}{p}\int_{\mathbb{S}^{n-1}}f(v)^pdS_{p, \gamma_n}(K, v),
 \end{eqnarray*}
where $h_t(v)=(h_K^p(v)+tf^p(v))^{\frac{1}{p}}$ for $p\neq 0$ and a continuous function $f: \mathbb{S}^{n-1}\rightarrow \mathbb{R}$, and $[h_t]$ is the Wulff shape, see \eqref{2.4}, generated by $h_t$. It was shown in \cite{73} that the $L_p$ Gaussian surface area measure of a convex body $K\in\mathcal{K}_o^n$ has the following integral representation 
\begin{eqnarray}\label{1.1a}
	S_{p, \gamma_n}(K, \eta)=\frac{1}{\left( \sqrt{2\pi}\right) ^n}\int_{\nu_K^{-1}(\eta)}e^{-\frac{|x|^2}{2}}\langle x, \nu_K(x)\rangle^{1-p}d\mathcal{H}^{n-1}(x)
\end{eqnarray}
for each Borel measure $\eta\subset \mathbb{S}^{n-1}$, where $\nu_K$ is the Gauss map and $\nu_K^{-1}$ is the inverse Gauss map (see Section 2 for the details of their definitions). Obviously, $S_{1, \gamma_n}(K, \cdot)$ is just the classical Gaussian surface area measure $S_{\gamma_n}(K, \cdot)$ introduced in \cite{71}. The $L_p$ Gaussian Minkowski problem can be stated in the following way:

 \noindent{\bf The $L_p$ Gaussian Minkowski problem.} Given a finite Borel measure $\mu$ on the unit sphere $\mathbb{S}^{n-1}$, what are the necessary and sufficient conditions on $\mu$ so that there exists a convex body $K\in\mathcal{K}_o^n$ such that
\begin{eqnarray*}
	S_{p ,\gamma_n}(K, \cdot)=\mu?
\end{eqnarray*} 
If $K$ exists, to what extent is it unique?

When $K$ is sufficiently smooth and the given measure $\mu$ has a density function $f$, the $L_p$ Gaussian Minkowski problem is equivalent to solving the following Monge-Amp\`{e}re type equation on $\mathbb{S}^{n-1}$:
 \begin{eqnarray}\label{1.1}
\frac{1}{(\sqrt{2\pi})^n}h^{1-p}_{K}e^{-\frac{\left|\nabla h_K\right|^2+h^2_K }{2}}\det\left(\nabla^2h_K+h_KI\right)=f,
\end{eqnarray}
where $h_K:\mathbb{S}^{n-1}\rightarrow \mathbb{R}$ is the support function of $K$, $\nabla$ and $\nabla^2$ are gradient and Hessian operators with respect to an orthonormal frame on $\mathbb{S}^{n-1}$, and $I$ is the identity matrix.

In \cite{73}, the existence of solutions to the $L_p$ Gaussian Minkowski problem for $p>0$ was established below.

{\it \noindent{\bf Theorem A \cite{73}.}~For $p> 0$, let $\mu$ be a non-zero finite Borel measure on $\mathbb{S}^{n-1}$ and not concentrate in any closed hemisphere. Then there exists a convex body $K\in \mathcal{K}^n_o$ such that
	\begin{eqnarray*}
	S_{p, \gamma_n}(K, \cdot)=\frac{S_{p, \gamma_n}(K, \mathbb{S}^{n-1})}{\mu(\mathbb{S}^{n-1})}\mu.
\end{eqnarray*}}

Since the $L_p$ Gaussian surface area measure does not have any homogeneity, the normalizing factor $S_{p, \gamma_n}(K, \mathbb{S}^{n-1})/\mu(\mathbb{S}^{n-1})$ in the equation above can not be removed in general. The same phenomenon occurs in solutions to the Orlicz Minkowski problems; see for example \cite{50, 51, 52, 56, 60, 61}. These solutions adjoint constant factors are often referred to as the normalization solutions for this type of Minkowski problems involving non-homogeneous geometry and probability measures.

The first objective of this paper is to consider the existence of symmetric solutions to the $L_p$ Gaussian Minkowski problem for $p\leq 0$. Namely, we obtain the following normalization solution in the origin-symmetric case.

{\it \noindent{\bf Theorem 1.1.}~Suppose $p\leq 0$ and $\mu$ is a non-zero finite even Borel measure on $\mathbb{S}^{n-1}$ that is not concentrated on any great subspheres. Then there exists an origin-symmetric convex body $K\in \mathcal{K}_e^n$ with $\gamma_n(K)=\frac{1}{2}$ such that
	\begin{eqnarray*}
	S_{p, \gamma_n}(K, \cdot)=\frac{S_{p, \gamma_n}(K, \mathbb{S}^{n-1})}{\mu(\mathbb{S}^{n-1})}\mu.
\end{eqnarray*}} 

The methods of proving the case $p<0$ in Theorem 1.1 are the variational method and an approximation argument. We first prove the existence result of solution in Theorem 1.1 when the given measure $\mu$ vanishes on all great subspheres of $\mathbb{S}^{n-1}$. In particular, this means that when the measure $\mu$ has a density function $f$, the eqaution \eqref{1.1} has a weak solution. Then it follows from the standard regularity theory of Monge-Amp\`{e}re equation that the equation \eqref{1.1} for $p<0$ has a smooth solution in the origin-symmetric case. Finally, an approximation argument is used to obtain the existence of the solution that we want.

It is well-know that the $L_0$ Minkowski problem is the logarithmic Minkowski problem the existence and uniqueness  of which are extremely difficult and remain open, see e.g., \cite{9, 10, 11, 16, 17, 19, 20, 21, 23, 32}. Obviously, it is of particularly meaningful to study the $L_0$ Gaussian Minkowski problem. As in \cite{9}, we prove the $L_0$ Gaussian Minkowski problem admits an origin-symmetric weak solution, see Theorem 1.1, by studying the relevant Monge-Amp\`{e}re type functional and using an approximation method.

Due to the lack of homogeneity of the $L_p$ Gaussian surface area measure, it is very difficult to remove the normalizing factor $S_{p, \gamma_n}(K, \mathbb{S}^{n-1})/\mu(\mathbb{S}^{n-1})$ in Theorem A. In \cite{73}, a degree-theoretic approach together with an appoximating argument was used to overcome this difficulty in the even case. Namely, a non-normalized solution (without the normalizing factor) for the even $L_p$ Gaussian Minkowski problem where $p\geq 1$ was established in \cite{73} as follows.

{\it \noindent{\bf Theorem B \cite{73}.}~For $p\geq 1$, let $\mu$ be a non-zero finite even Borel measure on $\mathbb{S}^{n-1}$ and not concentrate in any closed hemisphere with $|\mu|<\sqrt{\frac{2}{\pi}}r^{-p}ae^{-\frac{a^2}{2}}$ where $r$ and $a$ are chosen such that $\gamma_n(rB)=\gamma_n(P)=\frac{1}{2}$. Here $P$ is the symmetric strip $\{x\in \mathbb{R}^n: |x_1|\leq a\}$ and $B$ is the unit ball in $\mathbb{R}^n$. Then there exists a unique convex body $K\in \mathcal{K}_e^n$ with $\gamma_n(K)>\frac{1}{2}$ such that 
	\begin{eqnarray*}
	S_{p, \gamma_n}(K, \cdot)=\mu.
\end{eqnarray*}} 

The second aim of this paper is to give the following existence result for Theorem B without symmetry hypothesis on the given measure $\mu$. 

{\it \noindent{\bf Theorem 1.2.}~Let $\mu$ be a non-zero  finite Borel measure on $\mathbb{S}^{n-1}$ that is not concentrated in any closed hemisphere and $|\mu|<\left( \frac{n}{2}\right)^{1-p} \left( \frac{1}{\sqrt{2\pi}}\right)^p $ for $p\geq 1$. Then there exists a  convex body $K\in \mathcal{K}_o^n$ with $\gamma_n(K)>\frac{1}{2}$ such that 
	\begin{eqnarray*}
		S_{p, \gamma_n}(K, \cdot)=\mu.
\end{eqnarray*}}   

We will use the degree-theoretic method to establish Theorem 1.2, which is motivated by the works \cite{71, 73}. The crucial ingredient of using the degree theory is that a degree related to the equation \eqref{1.1} is well defined, where the $L_p$ Gaussian isoperimetric inequality proved in the following Theorem 5.1 plays an important role. Once the degree is well-defined, it follows from the fact that the degree remains invariant under continuous deformations of the equation \eqref{1.1} that the smooth solution of \eqref{1.1} exists by taking advantage of the existence of constant solutions of the equation \eqref{1.1} with $f$ being very small constant. Then, an approximation argument is used to give the weak solution in Theorem 1.2.

Finally, we also concentrate on the existence and uniqueness of asymmetric solutions to the $L_p$ Gaussian Minkowski problem in the smooth category when $p>n$. The existence is shown using the continuity method in PDEs, yet the uniqueness is obtained by the strong maximum principle in the uniformly elliptic equations. These results are stated as follows.

{\it \noindent{\bf Theorem 1.3.}~Let $f$ be a positive smooth function on $\mathbb{S}^{n-1}$. Then for $p> n$, there exists a unique smooth convex body $K$ solving the equation \eqref{1.1}.}

\ \ From a trivial calculation, it follows that the density of $S_{p, \gamma_n}(B_r, \cdot)$ of a centered ball of radius $r$ is $\frac{1}{(\sqrt{2\pi})^n}e^{-\frac{r^2}{2}}r^{n-p}$. Because of the behavior of $e^{-\frac{r^2}{2}}r^{n-p}\rightarrow 0 $ in the case $p<n$ as $r\rightarrow 0$ and $r\rightarrow \infty$, the uniqueness of solutions to the $L_p$ Gaussian Minkowski problem for $p<n$ does not hold without other conditions. Therefore, when $\gamma_n(K), \gamma_n(L)\geq\frac{1}{2}$ for $K, L\in \mathcal{K}_o^n$, the uniqueness for $p\geq 1$ was established in \cite{73}. Nevertheless, $e^{-\frac{r^2}{2}}r^{n-p}$ is monotonically decreasing on $r$ and $e^{-\frac{r^2}{2}}r^{n-p}\rightarrow 0$ (resp. $\infty$) as $r\rightarrow \infty$ (resp. $0$) when $p> n$. This means that solutions to the $L_p$ Gaussian Minkowski problem for $p> n$ are also unique if there are no additional condition that the Gaussian volume of convex bodies is greater than or equal to $1/2$. Therefore, we will use the strong maximum principle to obtain the uniqueness result in Theorem 1.3. Moreover, compared to the existence of smooth solution in Theorem 5.2, we also see that there is not the constraint of $f$ with $|f|_{L_1}<\left( \frac{n}{2}\right)^{1-p} \left( \frac{1}{\sqrt{2\pi}}\right)^p$ in the conditions of Theorem 1.3. This may be an advantage of the continuity method.

 The rest of this paper is organized as follows. In Section 2, we recall some notations and basics. In Section 3, we discuss the regularity of solutions to the $L_p$ Gaussian Minkowski problem based on the regularity results by Caffarelli \cite{5, 75}. The proof of Theorem 1.1 will be completed in Section 4. Section 5 is devoted to prove Theorem 1.2 whereas Section 6 is devoted to prove Theorem 1.3.

 \section{\bf Preliminaries}
\indent
In this section, we recall some notions and results which will be used later. More detailed information can be found in \cite{53, 64, 73, 76, 77}.

For a convex body $K\in \mathcal{K}_o^n$, the support function $h=h_K: \mathbb{R}^n\rightarrow \mathbb{R}$ and the radial function $\rho=\rho_K: \mathbb{R}^n\setminus\{o\}\rightarrow \mathbb{R}$ are respectively defined by
\begin{eqnarray*}
	h(x)=\max\{\langle x, y\rangle: y\in K\},~~ \rho(y)=\max\{\lambda: \lambda y\in K\}.
\end{eqnarray*}

The polar body of $K\in \mathcal{K}_o^n$ is defined by
\begin{eqnarray}\label{2.1z}
	K^*=\left\lbrace x\in \mathbb{R}^n: \langle x, y \rangle\leq 1~\mbox{for all}~y\in K\right\rbrace.
\end{eqnarray}
From \eqref{2.1z}, it easily follows that $K^*\in \mathcal{K}_o^n$ and $(K^*)^*=K$. In addition,
\begin{eqnarray}\label{2.2z}
	\rho_K=1/h_{K^*}, ~~~~~h_K=1/\rho_{K^*}.
\end{eqnarray}

For $K\in \mathcal{K}_o^n$ with $h_K(v_{\max})=\max_{v\in\mathbb{S}^{n-1}}h_K(v)$ and $\rho_K(u_{\min})=\min_{u\in\mathbb{S}^{n-1}}\rho_K(u)$, these facts found in \cite{64} that
\begin{eqnarray}\label{2.1}
	\max_{\mathbb{S}^{n-1}}h_K=\max_{\mathbb{S}^{n-1}}\rho_K,~~\min_{\mathbb{S}^{n-1}}h_K=\min_{\mathbb{S}^{n-1}}\rho_K, 
\end{eqnarray}
\begin{eqnarray}\label{2.2}
	h_K(v)\geq \langle v, v_{\max}\rangle h_K(v_{\max}) ~~\mbox{for any}~~v\in \mathbb{S}^{n-1},
\end{eqnarray}
\begin{eqnarray}\label{2.3}
	\rho_K(u)\langle u, u_{\min}\rangle\leq \rho_K(u_{\min}) ~~\mbox{for any}~~u\in \mathbb{S}^{n-1}
\end{eqnarray}
are of critical importance.

For a sequence of convex bodies $K_i \in\mathcal{K}_o^n$, we say that $K_i\rightarrow K_0\in \mathcal{K}_o^n$ with respect to the Hausdorff metric, provided when $i\rightarrow \infty$
\begin{eqnarray*}
	\|h_{K_i}-h_{K_0}\|_\infty=\max_{v\in \mathbb{S}^{n-1}}|h_{K_i}(v)-h_{K_0}(v)|\rightarrow 0,
\end{eqnarray*}
or equivalently
\begin{eqnarray*}
	\|\rho_{K_i}-\rho_{K_0}\|_\infty=\max_{u\in \mathbb{S}^{n-1}}|\rho_{K_i}(u)-\rho_{K_0}(u)|\rightarrow 0.
\end{eqnarray*}

The set of continuous functions on the unit sphere $\mathbb{S}^{n-1}$ will be denoted by $C(\mathbb{S}^{n-1})$, the subset of positive functions from $C(\mathbb{S}^{n-1})$ by $C^+(\mathbb{S}^{n-1})$, and the subset of even functions from $C^+(\mathbb{S}^{n-1})$ by $C_e^+(\mathbb{S}^{n-1})$.

For each $f\in C^+(\mathbb{S}^{n-1})$, the Wulff shape $\left[ f\right] $ generated by $f$ is a convex body defined by
\begin{eqnarray}\label{2.4}
	\left[ f\right]=\bigcap_{v\in \mathbb{S}^{n-1}}\left\lbrace x\in \mathbb{R}^{n}: \langle v, x\rangle\leq f(v)\right\rbrace. 
\end{eqnarray}
It is apparent that $h_{\left[ f\right]} \leq f $ and $\left[ h_K\right]=K$ for each $K\in \mathcal{K}_o^n$. 

For each $f\in C^+(\mathbb{S}^{n-1})$, the convex hull $\langle f\rangle$ generated by $f$ is also a convex body defined by 
\begin{eqnarray}\label{2.7z}
\langle f \rangle=\mbox{conv}\left\lbrace  f(u)u: u\in \mathbb{S}^{n-1}\right\rbrace.
\end{eqnarray}
Obviously, $\rho_{\langle f \rangle}\geq f$ and $\langle \rho_K\rangle=K$ for each $K\in \mathcal{K}_o^n$.

The Wulff shape $[f]$ of a continuous function $f\in C^+(\mathbb{S}^{n-1})$ and the convex hull $\langle 1/f\rangle$ generated by its reciprocal $1/ f$ are related by, see \cite[Lemma 2.8]{53},
\begin{eqnarray}\label{2.8z}
	[f]^*=\langle 1/f \rangle.
\end{eqnarray}

Let $K, L\in \mathcal{K}_o^n$ and $a, b\geq 0$. The $L_p$ Minkowski combination, $a\cdot K+_p b\cdot L\in \mathcal{K}_o^n$, is defined by the Wulff shape
\begin{eqnarray*}
a\cdot K+_pb\cdot L=\bigcap_{v\in \mathbb{S}^{n-1}}\left\lbrace x\in \mathbb{R}^n: \langle x, v \rangle\leq \left(ah_K(v)^p+bh_L(v)^p\right)^{\frac{1}{p}}\right\rbrace 
\end{eqnarray*}
for $p\neq 0$ and
\begin{eqnarray*}
	a\cdot K+_0b\cdot L=\bigcap_{v\in \mathbb{S}^{n-1}}\left\lbrace x\in \mathbb{R}^n: \langle x, v \rangle\leq h_K(v)^ah_L(v)^b\right\rbrace
\end{eqnarray*}
for $p=0$.

The supporting hyperplane $H_K(v)$ of $K\in\mathcal{K}_o^n$ for each $v\in \mathbb{S}^{n-1}$ is defined by
\begin{eqnarray*}
	H_K(v)=\{x\in \mathbb{R}^n: \langle x, v\rangle=h_K(v)\}.
\end{eqnarray*}

Let $\partial K$ denote the boundary of a convex body $K$. For $\sigma\subset \partial K$ with $K\in\mathcal{K}_o^n$, the spherical image of $\sigma$ is defined by
\begin{eqnarray*}
	\boldsymbol{\nu}_K(\sigma)=\left\{v\in \mathbb{S}^{n-1}: x\in H_K(v) ~\mbox{for some}~ x\in \sigma\right\}\subset \mathbb{S}^{n-1}.
\end{eqnarray*} 
Let $\sigma_K\subset \partial K$ be set consisting of all $x\in \partial K$ for which the set $\boldsymbol{\nu}_K(\{x\})$, often abbreviated as $\boldsymbol{\nu}_K(x)$, contains more than a single element. It is well known that $\mathcal{H}^{n-1}(\sigma_K)=0$ (see Schneider \cite [p.84]{77}). Let ${\nu}_K(x)$ be the unique element in $\boldsymbol{\nu}_K(x)$ for each $x\in \partial K\setminus\sigma_K$. Then we can define the function
\begin{eqnarray*}
	\nu_K: \partial K\setminus \sigma_K\rightarrow \mathbb{S}^{n-1},
\end{eqnarray*}
which is called the spherical image map (also known as the Gauss map) of $K$. It will occasionally be convenient to abbreviate $\partial K\setminus \sigma_K$ by $\partial'K$.

For $\eta \subset\mathbb{S}^{n-1}$, the reverse spherical image $\boldsymbol{\nu}^{-1}_K$, of $K\in\mathcal{K}_o^n$ at $\eta$, is defined by
\begin{eqnarray*}
	\boldsymbol{\nu}^{-1}_K(\eta)=\left\{x\in \partial K: x\in H_K(v) ~\mbox{for some}~ v\in \eta\right\}\subset \partial K. 
\end{eqnarray*} 
The set $\eta_K\subset \mathbb{S}^{n-1}$ consisting of all $v\in \mathbb{S}^{n-1}$ for which the set $\boldsymbol{\nu}^{-1}_K(v)=\boldsymbol{\nu}^{-1}_K(\{v\})$ contains more than a single element is of $\mathcal{H}^{n-1}$-measure $0$ (see Schneider \cite[Theorem 2.2.11]{77}). For $v\in \mathbb{S}^{n-1}\setminus\eta_K$, $\boldsymbol{\nu}^{-1}_K(v)$ has the unique element denoted by ${\nu}^{-1}_K(v)$. Hence we define the reverse spherical image map
\begin{eqnarray*}
	{\nu}^{-1}_K: \mathbb{S}^{n-1}\setminus\eta_K\rightarrow \partial K,
\end{eqnarray*} 
which is continuous from Lemma 2.2.12 of Schneider \cite{77}. 

Let $\omega\subset \mathbb{S}^{n-1}$ be a Borel set. The radial Gauss image of $K\in\mathcal{K}_o^n$ at $\omega$, denoted by $\boldsymbol{\alpha}_K(\omega)$, is defined by
\begin{eqnarray*}
	\boldsymbol{\alpha}_K(\omega)=\left\{v\in\mathbb{S}^{n-1}: \rho_K(u)u\in H_K(v) ~\mbox{for some}~ u\in \omega\right\}\subset \mathbb{S}^{n-1}. 
\end{eqnarray*}  
When $\omega=\{u\}$ is a singleton, we usually write $\boldsymbol{\alpha}_K(u)$ instead of $\boldsymbol{\alpha}_K(\{u\})$. Let $\omega_K=\{u\in \mathbb{S}^{n-1}: \rho_K(u)u\in\sigma_K\}\subset \mathbb{S}^{n-1}$. Clearly, $\boldsymbol{\alpha}_K(u)$ contains more than one element for each $u\in \omega_K$. Note that $\mathcal{H}^{n-1}(\omega_K)=0$ from Theorem 2.2.5 of \cite {77}. The radial Gauss map of $K$, denoted by $\alpha_K$, is the map defined on $\mathbb{S}^{n-1}\setminus\omega_K$ that takes each point $u$ in its domain to the unique vector in $\boldsymbol{\alpha}_K(u)$. Hence $\alpha_K$ is defined almost everywhere on $\mathbb{S}^{n-1}$ with respect to the spherical Lebesgue measure.

Let $\eta \subset\mathbb{S}^{n-1}$ be a Borel set. The reverse radial Gauss image $\boldsymbol{\alpha}^\ast_K(\eta)$, of $K\in\mathcal{K}_o^n$ at $\eta$, is defined by
\begin{eqnarray*}
	\boldsymbol{\alpha}^\ast_K(\eta)=\left\{u\in\mathbb{S}^{n-1}: \rho_K(u)u\in H_K(v) ~\mbox{for some}~ v\in \eta\right\}\subset \mathbb{S}^{n-1}. 
\end{eqnarray*} 
Similarly, we write $\boldsymbol{\alpha}^\ast_K(v)$ instead of $\boldsymbol{\alpha}^\ast_K(\{v\})$ for $\eta=\{v\}$. Obviously,  For $v\in \mathbb{S}^{n-1}\setminus\eta_K$, $\boldsymbol{\alpha}^{\ast}_K(v)$ has the unique element denoted by ${\alpha}^{\ast}_K(v)$. Thus we define the reverse radial Gauss image map
\begin{eqnarray*}
	{\alpha}^{\ast}_K: \mathbb{S}^{n-1}\setminus\eta_K\rightarrow \mathbb{S}^{n-1}.
\end{eqnarray*} 
Therefore, $\alpha^{\ast}_K$ is defined almost everywhere on $\mathbb{S}^{n-1}$ since the set $\eta_K$ has spherical 
Lebesgue measure $0$.

The surface area measure $S(K, \cdot)$ of a convex body $K$ is a Borel measure on  $\mathbb{S}^{n-1}$ defined, for a Borel set $\omega\subset \mathbb{S}^{n-1}$, by
\begin{eqnarray*}
	S(K,\omega)=\mathcal{H}^{n-1}(\nu_K^{-1}(\omega))=\mathcal{H}^{n-1}(\{x\in \partial K: \nu_K(x)\cap \omega\neq \emptyset\}),
\end{eqnarray*}
which is equivalently defined by

\begin{eqnarray}\label{2.9z}
	\int_{\mathbb{S}^{n-1}}f(v)dS(K, v)=\int_{\partial^\prime K}f(\nu_K(x))d\mathcal{H}^{n-1}(x),
\end{eqnarray}
for every $f\in C(\mathbb{S}^{n-1})$.

For each bounded Lebesgue integrable function $f:\mathbb{S}^{n-1}\rightarrow \mathbb{R}$ and $K\in\mathcal{K}_o^n$, the following fact was established in \cite{53}:
\begin{eqnarray}\label{2.10z}
	\int_{\mathbb{S}^{n-1}}f(u)\rho_K(u)^ndu=\int_{\partial^\prime K}\langle x, \nu_K(x)\rangle f\left(\frac{x}{|x|}\right)d\mathcal{H}^{n-1}(x).
\end{eqnarray}

Let $f \in C(\mathbb{S}^{n-1})$ and $\varrho_0, \varrho_t\in C^+(\mathbb{S}^{n-1})$. For each $t\in (-\delta, \delta)$ with $\delta>0$, define
\begin{eqnarray}\label{2.11z}
\log\varrho_{t}(u)=\log\varrho_0(u)+tf(u)+o(t, u),
\end{eqnarray}
where $o(t, \cdot): \mathbb{S}^{n-1}\rightarrow \mathbb{R}$ is continuous and $\lim_{t\rightarrow 0}\frac{o(t, \cdot)}{t}=0$ uniformly on $\mathbb{S}^{n-1}$.

The following facts from \cite[Lemmas 4.1-4.2]{53} will be required.

{\it \noindent{\bf Lemma 2.1 \cite{53}.}~Let $f \in C(\mathbb{S}^{n-1})$, $\varrho_0 \in C^+(\mathbb{S}^{n-1})$ and $\varrho_t$ be given as in \eqref{2.11z}. Then
	
 {\bf (I)}
	\begin{eqnarray}\label{2.12z}
	\lim_{t\rightarrow 0}\frac{\log h_{\langle \varrho_t \rangle}(v)-\log h_{\langle \varrho_0 \rangle}(v)}{t}=f(\alpha^*_{\langle \varrho_0 \rangle}(v))
\end{eqnarray} 
holds for almost all $v\in \mathbb{S}^{n-1}$. Moreover, there exist $\delta_0>0$ and $M>0$ so that
\begin{eqnarray}\label{2.13z}
|\log h_{\langle \varrho_t \rangle}(v)-\log h_{\langle \varrho_0 \rangle}(v)|\leq M|t|,
\end{eqnarray} 	
for all $v\in \mathbb{S}^{n-1}$ and all $t\in (-\delta_0, \delta_0)$.	
	
 {\bf (II)}	
	\begin{eqnarray}\label{2.14z}
	\lim_{t\rightarrow 0}\frac{ h^{-1}_{\langle \varrho_t \rangle}(v)-h^{-1}_{\langle \varrho_0 \rangle}(v)}{t}=-h^{-1}_{\langle \varrho_0 \rangle}(v)f(\alpha^*_{\langle \varrho_0 \rangle}(v))
\end{eqnarray} 
holds for almost all $v\in \mathbb{S}^{n-1}$. Moreover, there exist $\delta_0>0$ and $M>0$ such that
\begin{eqnarray}\label{2.15z}
	|h^{-1}_{\langle \varrho_t \rangle}(v)-h^{-1}_{\langle \varrho_0 \rangle}(v)|\leq M|t|,
\end{eqnarray} 	
for all $v\in \mathbb{S}^{n-1}$ and all $t\in (-\delta_0, \delta_0)$. }
 
The following is the variational formula for Gaussian volume that shows that the $L_p$ Gaussian surface area measure is differentials of the Gaussian volume, which was discovered in \cite{73}.
 
{\it \noindent{\bf Lemma 2.2 \cite{73}.}~Let $p\neq 0$, $K\in\mathcal{K}_o^n$ and $f\in C(\mathbb{S}^{n-1})$. For sufficiently small $\delta>0$ and each $t\in (-\delta, \delta)$, define $h_t\in C^+(\mathbb{S}^{n-1})$ by $h_t=(h_K^p+tf^p)^{\frac{1}{p}}$. Then
\begin{eqnarray*}
\lim_{t\rightarrow 0}\frac{\gamma_n([h_t])-\gamma_n(K)}{t}=\frac{1}{p}\int_{\mathbb{S}^{n-1}}f(v)^pdS_{p, \gamma_n}(K, v).
\end{eqnarray*} 		} 

{\it \noindent{\bf Lemma 2.3 \cite{73}.}~Let $p\geq 1$ and  $K, L\in \mathcal{K}_o^n$ with $\gamma_n(K), \gamma_n(L)\geq\frac{1}{2}$. Then 
	\begin{eqnarray*}
		S_{p, \gamma_n}(K, \cdot)=S_{p, \gamma_n}(L, \cdot) \Longrightarrow K=L.
	\end{eqnarray*}} 
  
The $L_p$ Gaussian surface area measure is weakly convergent with respect to Hausdorff metric (see \cite{73}): If $K_i \in\mathcal{K}_o^n$ such that $K_i\rightarrow K_0\in \mathcal{K}_o^n$ as $i\rightarrow\infty$, then $S_{p, \gamma_n}(K_i, \cdot) \rightarrow S_{p, \gamma_n}(K_0, \cdot)$ weakly.

\section{\bf Regularity}
\indent In this section, we aim to study the regularity of solutions to the $L_p$ Gaussian Minkowski problem, which is used to prove the existence of smooth solutions to the problem. Let us first collect some basic definitions and facts about regularity theory. We can consult \cite{78} and \cite{79} for some details.

  Let $K$ be a convex body. If $\partial K$ contains no segment then we say that $K$ is strictly convex; if $K$ has a unique tangential hyperplane at $x\in \partial K$ then we say that $x$ is a $C^1$-smooth point. Clearly, $h_K$ is $C^1$ on $\mathbb{S}^{n-1}$ if and only if $K$ is strictly convex. Moreover, $\partial K$ is $C^1$ if and only if each $x\in \partial K$ is $C^1$-smooth.
  
  Let $\Omega$ be a convex set in $\mathbb{R}^{n}$. If $y=\lambda x_1+(1-\lambda)x_2$ for $x_1, x_2\in \Omega$ and $\lambda\in (0, 1)$ implies $x_1=x_2=y$, then we say that $y\in \Omega$ is an extremal point. Note that if $\Omega$ is compact and convex, then $\Omega$ is the convex hull of its extremal points.
  
  The normal cone of a convex body $K$ at $z\in K$ is defined by
  \begin{eqnarray*}
  	N(K, z)=\left\lbrace x\in \mathbb{R}^{n}: \langle x, y\rangle\leq \langle x, z\rangle~\mbox{for all}~y\in K\right\rbrace, 
  \end{eqnarray*}
  or, equivalently by
  \begin{eqnarray*}
  	N(K, z)=\left\lbrace x\in \mathbb{R}^{n}:h_K(x)=\langle x, z \rangle\right\rbrace.
  \end{eqnarray*}
  If $z\in \mbox{int} K$ then we see $N(K, z)=\left\lbrace  o \right\rbrace$, and if $z\in \partial K$ then $\mbox{dim}N(K, z)\geq 1$.
  
  The face of a convex body $K$ with outer normal $x\in \mathbb{R}^{n}$ is defined as 
  \begin{eqnarray}\label{3.1}
  	F(K, x)=\left\lbrace y\in K:h_K(x)=\langle x, y \rangle\right\rbrace,
  \end{eqnarray}
  which lies in $\partial K$ provided $x\neq o$, and 
  \begin{eqnarray}\label{3.2}
  	F(K, x)=\partial h_K(x).
  \end{eqnarray}
  Here $\partial h_K(x)$ is the subgradient of $h_K$, namely
  \begin{eqnarray*}
  	\partial h_K(x)=\left\lbrace z\in \mathbb{R}^{n}: h_K(y)\geq h_K(x)+\langle z, y-x\rangle~\mbox{for each } y\in K\right\rbrace, 
  \end{eqnarray*}
  which is a non-empty compact convex set. Note that $h_K(x)$ is differentiable at $x$ if and only if $\partial h_K(x)$ consists of exactly one vector which is the gradient, denoted by $Dh_K(x)$, of $h_K$ at $x$.
  
 When $h_K$ is viewed as restricted to the unit sphere $\mathbb{S}^{n-1}$, the gradient of $h_K$ on $\mathbb{S}^{n-1}$ is denoted by $\nabla h_K$. Since $h_K$ is differentiable at $\mathcal{H}^n$ almost all points in $\mathbb{R}^n$ and is positively homogeneous of degree $1$, $h_K$ is differentiable for $\mathcal{H}^{n-1}$ almost all points of $\mathbb{S}^{n-1}$. Let $h_K$ be differentiable at $v\in \mathbb{S}^{n-1} $, where $v=\nu_K(x)$ is an outer unit normal vector at $x\in \partial K$. Then we have
\begin{eqnarray}\label{3.3}
	x=Dh_K(v)=\nu_K^{-1}(v).
\end{eqnarray}
This implies
\begin{eqnarray}\label{3.4}
	h_K(v)=h_K(\nu_K(x))=\langle x, \nu_K(x)\rangle=\langle Dh_K(v), v\rangle,
\end{eqnarray}
\begin{eqnarray}\label{3.5}
	Dh_K(v)=\nabla h_K(v)+h_K(v)v,
\end{eqnarray}
\begin{eqnarray}\label{3.6}
	|Dh_K(v)|=\sqrt{h_K(v)^2+|\nabla h_K(v)|^2}.
\end{eqnarray} 
 
The following is to recall the notions and facts about Monge-Amp\`{e}re measure from the survey by Trudinger and Wang \cite{80}.

  Given a convex function $\varphi$ defined in an open convex set $\Omega$ of $\mathbb{R}^{n}$, $D\varphi$ and $D^2\varphi$ denote its gradient and its Hessian, respectively. For any Borel subset $\vartheta \subset\Omega$, define 
  \begin{eqnarray*}
  	N_\varphi(\vartheta)=\bigcup_{x\in\vartheta}\partial \varphi(x).
  \end{eqnarray*}
  The Monge-Amp\`{e}re measure $\mu_\varphi$ is $\mu_\varphi(\vartheta)=\mathcal{H}^n\left( N_\varphi(\vartheta)\right)$. If the function $\varphi$ is $C^2$ smooth, then the subgradient $\partial \varphi$ coincide with the gradient $D \varphi$. Thus
  \begin{eqnarray}\label{3.7}
  	\mu_\varphi(\vartheta)=\mathcal{H}^n\left( D \varphi(\vartheta)\right)=\int_\vartheta \mbox{det}\left(D^2 \varphi\right)d\mathcal{H}^n.
  \end{eqnarray}
  Note that the surface area measure $S(K, \cdot)$ of a convex body $K$ in $\mathbb{R}^{n}$ is a Monge-Amp\`{e}re type measure with $h_K$ restricted to the unit sphere $\mathbb{S}^{n-1}$ since it satisfies
  \begin{eqnarray}\label{3.8}
  	S(K, \omega)=\mathcal{H}^{n-1}\left(\bigcup_{v\in\omega}F(K, v)\right)=\mathcal{H}^{n-1}\left( \bigcup_{v\in\omega} \partial h_K(v)\right)=\mu_{h_K}(\omega)
  \end{eqnarray}
  for any Borel $\omega\subset \mathbb{S}^{n-1}$.
  
  We say that a convex function $\varphi$ is the solution of a Monge-Amp\`{e}re equation in the sense of measure (or in the Aleksandrov sense), if it solves the corresponding integral formula for $\mu_\varphi$ instead of the original formula for $\mbox{det}(D^2 \varphi)$.
  
  The first step toward studying the regularity of solutions to the $L_p$ Gaussian Minkowski problem is to convert the original Monge-Amp\`{e}re type equation (\ref{1.1}) on the unit sphere $\mathbb{S}^{n-1}$ into a Euclidean Monge-Amp\`{e}re type equation on $\mathbb{R}^{n-1}$. In view of this, we consider the restriction of a solution $h$ of (\ref{1.1}) to the hyperplane tangential to $\mathbb{S}^{n-1}$ at $e\in \mathbb{S}^{n-1}$.

  {\it \noindent{\bf Lemma 3.1.}~Let $e\in \mathbb{S}^{n-1}, K\in \mathcal{K}_o^n$, and $\varphi: e^\perp\rightarrow \mathbb{R}$ with $\varphi(y)=h_K(y+e)$. If $h=h_K$ is a solution of (\ref{1.1}) for non-negative functions $f$, then $\varphi$ satisfies the standard Monge-Amp\`{e}re equation on $e^\perp$ in the sense of measure:
  	\begin{eqnarray}\label{3.9}
  		{\rm det} D^2\varphi(y)=(\sqrt{2\pi})^n\varphi(y)^{p-1}e^{\frac{\left| D\varphi(y)+(\varphi(y)-\langle D\varphi(y), y\rangle)\cdot e\right|^2}{2}}G(y),
  	\end{eqnarray}	
 where 	
 \begin{eqnarray*}
 	G(y)= \left(1+|y|^2\right)^{-\frac{n+p}{2}}f\left(\dfrac{e+y}{\sqrt{1+|y|^2}}\right).
 \end{eqnarray*}	
}
  
  {\bf Proof.}~Suppose that for $K\in \mathcal{K}_o^n$, $h=h_K$
  solves the equation (\ref{1.1}). Since
   \begin{eqnarray}\label{3.10}
  dS(K, v)=\det\left(\nabla^2h_K(v)+h(v)I\right)dv
  \end{eqnarray}
for $v\in \mathbb{S}^{n-1}$,  it follows from (\ref{1.1})  and  (\ref{3.6}) that
  \begin{eqnarray}\label{3.11}
  	dS(K, v)=(\sqrt{2\pi})^nfh_K^{p-1}(v)e^{\frac{|Dh_K(v)|^2}{2}}dv.
  \end{eqnarray}	
  
  For $e\in \mathbb{S}^{n-1}$, let $H_e$ be the hyperplane in $\mathbb{R}^{n}$ which is tangential to $\mathbb{S}^{n-1}$ at $e$, and $e^\perp$ the orthogonal complement of $\{re: r\in \mathbb{R}\}$ in $\mathbb{R}^n$. For $y\in e^\perp$, we have $y=\sum_{i=1}^{n-1}y_ie_i$, where $\{e_1, \cdots, e_{n-1}\}$ is a basis of $e^\perp$. Denoted by $\pi: e^\perp\rightarrow \mathbb{S}^{n-1}$ the radial projection from $H_e=e+e^\perp$ to $\mathbb{S}^{n-1}$, which is defined by
  \begin{eqnarray*}
  	\pi(y)=\frac{y+e}{\sqrt{1+|y|^2}}.
  \end{eqnarray*}	
  By the fact that 
  \begin{eqnarray}\label{3.12}
  	\langle \pi(y), e\rangle=\frac{1}{\sqrt{1+|y|^2}},
  \end{eqnarray}	
  we can calculate that the determinant of the Jacobin of the mapping $y\mapsto v=\pi (y)$ is
  \begin{eqnarray}\label{3.13}
  	\left| \mbox{Jac}\pi\right| =\left| \frac{dv}{dy}\right| =\left( \dfrac{1}{\sqrt{1+|y|^2}} \right)^n.
  \end{eqnarray}
  Let $\varphi: e^\perp\rightarrow \mathbb{R}$ be the restriction of $h_K$ on $H_e$. That is,
  \begin{eqnarray}\label{3.14}
  	\varphi(y)=h_K(y+e)=\sqrt{1+|y|^2}h_K(\pi(y)).
  \end{eqnarray}
  Then it is not hard to see from (\ref{3.1}) and (\ref{3.2}) that
  \begin{eqnarray}\label{3.15}
  	\partial \varphi(y)= \partial h_K(y+e)\left| \right.e^\perp=F(K, y+e)\left|\right. e^\perp=F(K, \pi(y))\left| \right.e^\perp.
  \end{eqnarray}
  It follows from the homogeneity of degree $1$ and the differentiability of $h_K$ that
  \begin{eqnarray*}
  	Dh_K(y+e)=Dh_K(v),
  \end{eqnarray*}
  where $v=\pi (y)$. From this, we see
  \begin{eqnarray*}
  	D\varphi(y)=Dh_K(y+e)\left| \right.e^\perp=Dh_K(v)\left| \right.e^\perp.
  \end{eqnarray*}
  Therefore we can assume that
  \begin{eqnarray}\label{3.16}
  	Dh_K(v)=D\varphi (y)-t e
  \end{eqnarray}
  for some undetermined constant $t\in \mathbb{R}$. By (\ref{3.4}), we see
  \begin{eqnarray}\label{3.17}
  	h_K(v)=\langle Dh_K(v), v\rangle.
  \end{eqnarray}
  Note also that
  \begin{eqnarray}\label{3.18}
  	v=\pi(y)=\frac{y+e}{\sqrt{1+|y|^2}},~~h_K(v)=\dfrac{\varphi(y)}{\sqrt{1+|y|^2}}.
  \end{eqnarray}
  Substituting (\ref{3.16}) and (\ref{3.18}) into (\ref{3.17}), it follows that
  \begin{eqnarray}\label{3.19}
  	t=\langle D\varphi(y), y\rangle-\varphi(y).
  \end{eqnarray}
  Thus combining (\ref{3.16}) with (\ref{3.19}) we have
  \begin{eqnarray}\label{3.20}
  	Dh_K(v)=D\varphi(y)+\left(\varphi(y)-\langle D \varphi(y), y\rangle\right) \cdot e.
  \end{eqnarray}
  For a Borel set $\vartheta \subset e^\perp$, based on (\ref{3.8}) it is easy to check that
  \begin{eqnarray}\label{3.21}
  	\mathcal{H}^{n-1}\left(\bigcup_{v\in \pi(\vartheta)}\left( F(K, v)\left| e^\perp\right. \right) \right)=\int_{\pi(\vartheta)}\langle v, e\rangle dS(K, v).
  \end{eqnarray}
  Thus it follows from (\ref{3.15}), (\ref{3.4}), (\ref{3.11}), (\ref{3.20}), (\ref{3.13}) and (\ref{3.14}) that
  \begin{eqnarray}\label{3.21z}
  	&&\int_\vartheta \mbox{det}D^2 \varphi (y)d\mathcal{H}^{n-1}(y) \nonumber\\
  	&=&\int_{\pi(\vartheta)}\langle v, e\rangle dS(K, v)\\
  	&=&(\sqrt{2\pi})^n\int_{\pi(\vartheta)}\langle v, e\rangle h_K(v)^{p-1} e^{\frac{|Dh_K(v)|^2}{2}}f(v)dv \nonumber\\
  	&=&(\sqrt{2\pi})^n\int_\vartheta \varphi(y)^{p-1}e^{\frac{\left| D\varphi(y)+(\varphi(y)-\langle D\varphi(y), y\rangle)\cdot e\right|^2}{2}}G(y)d\mathcal{H}^{n-1}(y),\nonumber
  \end{eqnarray}
where 
 \begin{eqnarray*}
	G(y)= \left(1+|y|^2\right)^{-\frac{n+p}{2}}f\left(\dfrac{e+y}{\sqrt{1+|y|^2}}\right).
\end{eqnarray*}
This implies that $\varphi$ satisfies (\ref{3.9}) on $e^\perp$.
 \hfill $\square$
  
  The following two lemmas from Caffarelli \cite{5, 75}, see also \cite{79}, play crucial role in the discussion of the regularity of solutions to the $L_p$ Gaussian Minkowski problem.
  
  {\it \noindent{\bf Lemma 3.2}(Caffarelli\cite{75}).~Let $\lambda_2>\lambda_1>0$, and let $\varphi$ be a convex function on an open bounded convex set $\Omega\subset \mathbb{R}^n$ such that
  	\begin{eqnarray*}
  		\lambda_1\leq \det D^2 \varphi \leq \lambda_2
  	\end{eqnarray*}	
  	in the sense of measure.
  	
  	\noindent(i) If $\varphi$ is non-negative and $W=\{y\in \Omega: \varphi(y)=0\}$ is not a point, then $W$ has no 
  	
  	extremal point in $\Omega$.
  	
  	\noindent(ii) If $\varphi$ is strictly convex, then $\varphi$ is $C^1$.
  }
  
  {\it \noindent{\bf Lemma 3.3}(Caffarelli\cite{5}).~For real functions $\varphi$ and $f$ on an open bounded convex set $\Omega\subset \mathbb{R}^n$, let $\varphi$ be strictly convex, and let $f$ be positive and continuous such that 
  	\begin{eqnarray*}
  		\det D^2 \varphi =f
  	\end{eqnarray*}	
  	in the sense of measure.
  	
  	\noindent(i) Each $z\in \Omega$ has an open ball $B\subset \Omega$ around $z$ such that the restriction of $\varphi$ to $B$ 
  	
  	is in $C^{1, \alpha}(B)$ for any $\alpha\in (0, 1)$.

  	\noindent(ii) If $f$ is in $C^\alpha (\Omega)$ for some $\alpha\in (0, 1)$, then each $z\in \Omega$ has an open ball $B\subset \Omega$ 
  	
  	around $z$ such that the restriction of $\varphi$ to $B$ is in $C^{2, \alpha}(B)$.
  }
  
  With the help of the above lemmas, we prove the following results of regularity.
  
 {\it \noindent{\bf Theorem 3.1.}~Suppose that $d\mu=fd\mathcal{H}^{n-1}$ with $0<c_1\leq f \leq c_2$ on $\mathbb{S}^{n-1}$, and $K\in\mathcal{K}_o^n$ satisfies $dS_{p, \gamma_{n}}(K, \cdot)=fd\mathcal{H}^{n-1}$ on $\mathbb{S}^{n-1}$. Then
 	
 	\noindent(i) $\partial K$ is $C^1$ and strictly convex, and $h_K$ is $C^1$ on $\mathbb{R}^n\setminus\{o\}$;
 	
 	\noindent(ii) if $f$ is continuous, then the restriction of $h_K$ to $\mathbb{S}^{n-1}$ is in $C^{1, \alpha}$ for any $\alpha\in (0, 1)$;
 	
 	\noindent(iii) if $f \in C^\alpha(\mathbb{S}^{n-1})$ for $\alpha\in (0, 1)$, then $h_K$ is $C^{2, \alpha}$ on $\mathbb{S}^{n-1}$.}

  {\bf Proof.} Let $e\in \mathbb{S}^{n-1}$ and $0<\tau< 1$. We define
  \begin{eqnarray*}
  	\Phi(e, \tau)=\{v\in\mathbb{S}^{n-1}: \langle v, e\rangle> \tau\}.
  \end{eqnarray*}
  Note that $h_K$ is continuous on $\mathbb{S}^{n-1}$ for $K\in \mathcal{K}_o^n$. Hence there exist $0<\tau_1<1$ and $\delta>0$ such that
  \begin{eqnarray*}
  	h_K(v)\geq\delta, ~\mbox{for}~v\in \mbox{cl}\Phi(e, \tau_1).
  \end{eqnarray*}
  Here, both $\tau_1$ and $\delta$ depend on $e$ and $K$. Moreover, there exists $0<\epsilon<1$ depending on $e$ and $K$ such that if some $v\in \mbox{cl}\Phi(e, \tau_1)$ is the outer normal at $x\in\partial K$, then
  \begin{eqnarray}\label{3.22}
  	\epsilon \leq |x| \leq \frac{1}{\epsilon}.
  \end{eqnarray}
  Recall that $\pi$ is defined by
  \begin{eqnarray*}
  	\pi(y)=\frac{y+e}{\sqrt{1+|y|^2}}, ~~y\in e^\bot.
  \end{eqnarray*}	
  Define 
  \begin{eqnarray*}
  	\Psi_e=\pi^{-1}(\Phi (e, \tau_1)).
  \end{eqnarray*}	
  Let $\varphi: e^\perp\rightarrow \mathbb{R}$ satisfy the assumption of Lemma 3.1. Then for $y\in \Psi_e$, it follows that from (\ref{3.3}), (\ref{3.20}) and (\ref{3.22}) that
  \begin{eqnarray}\label{3.23}
  	\epsilon\leq |D\varphi(y)+\left( \varphi(y)-\langle D\varphi(y), y\rangle\right) \cdot e|\leq \frac{1}{\epsilon}.
  \end{eqnarray}	
  Note that for $y\in \mbox{cl}\Psi_e$,
  \begin{eqnarray*}
  	\varphi(y)=\sqrt{1+|y|^2}h_K\left(\dfrac{e+y}{\sqrt{1+|y|^2}} \right) \geq \delta.
  \end{eqnarray*}
  In addition, it is easy to see that $\varphi$ also has an upper bound depending on $e$ and $K$ for $y\in \mbox{cl}\Psi_e$. Since it is assumed that for positive constants $c_1$ and $c_2$,  
  \begin{eqnarray*}
  	0<c_1\leq f \leq c_2.
  \end{eqnarray*}
  Thus it follows from Lemma 3.1 and (\ref{3.23}) that there exists $\zeta\in (0, 1)$ depending on $e$ and $K$ such that for $y\in \Psi_e$,
  \begin{eqnarray}\label{3.24}
  	\zeta\leq \det D^2\varphi (y)\leq \frac{1}{\zeta}.
  \end{eqnarray}
  
  We first show that $\partial K$ is $C^1$ for $K\in \mathcal{K}_o^n$. Namely, $\mbox{dim} N(K, z)=1$ for any 
  $z\in \partial K$. If not, then there exists a point $z_0\in \partial K$ such that $\mbox{dim} N(K, z_0)\geq 2$. Assume $e\in N(K, z_0) \cap \mathbb{S}^{n-1}$. By the definition of support function, and noting $z_0\in \partial K$, it follows that for $y\in \Psi_e$,
  \begin{eqnarray*}
  	\varphi(y)=h_K(y+e)\geq \langle y+e, z_0\rangle.
  \end{eqnarray*}
  It is not hard to verify that 
  \begin{eqnarray*}
  	y\in W:=\pi^{-1}(N(K, z_0)\cap \Phi (e, \tau_1))\Longleftrightarrow \varphi(y)=\langle y+e, z_0\rangle~\mbox{for}~y \in \Psi_e.
  \end{eqnarray*}
  Let $\iota(y)=\langle y+e, z_0\rangle$. Then
  \begin{eqnarray*}
  	\varphi(y)-\iota(y)\left\{
  	\begin{aligned}
  		&=&0~~\mbox{for}~&y\in W&\\
  		&>&0~~\mbox{for}~&y\in \Psi_e\setminus W.&
  	\end{aligned}
  	\right.	
  \end{eqnarray*}
  By (\ref{3.24}), and noting that $\iota$ is the first degree polynomial, we get that for $y\in \Psi_e$,
  \begin{eqnarray*}
  	\zeta\leq \det D^2(\varphi (y)-\iota(y))\leq \frac{1}{\zeta}.
  \end{eqnarray*}
  Note that
  \begin{eqnarray*}
  	W=\pi^{-1}(N(K, z_0)\cap \Phi (e, \tau_1))=\{y\in \Psi_e: \varphi (y)-\iota(y)=0\}
  \end{eqnarray*}
  is not a point since $\mbox{dim} W\geq 1$. Since the origin $o$ is an extremal point of $W$ by the choice of $e$, this is a contradiction from (i) of Lemma 3.2. 
  
  Next, we prove that $\varphi$ is strictly convex on $\mbox{cl} \Psi_e$ for $e\in \mathbb{S}^{n-1}$. It is easy to see that $\Psi_e$ is a convex set in $\mathbb{R}^{n-1}$. For $0<\lambda<1$ and $y_1, y_2\in \Psi_e$ with $y_1\neq y_2$, we assume that $e+(\lambda y_1+(1-\lambda)y_2)$ is an outer normal at $z\in \partial K$. Namely,
  \begin{eqnarray*}
  	e+(\lambda y_1+(1-\lambda)y_2)\in N(K, z).
  \end{eqnarray*}
  Note that
  \begin{eqnarray*}
  	e+y_1\notin N(K, z)~~\mbox{and}~~e+y_2\notin N(K, z)
  \end{eqnarray*}
  since $z\in \partial K$ is a smooth point. Thus,
  \begin{eqnarray*}
  	\varphi(y_i)=h_K(e+y_i)>\langle z,e+y_i\rangle.
  \end{eqnarray*}
  From this, we deduce
  \begin{eqnarray*}
  	\lambda \varphi(y_1)+(1-\lambda) \varphi (y_2)&>&
  	\langle z, e+(\lambda y_1+(1-\lambda)y_2)\rangle\\
  	&=&h_K(e+(\lambda y_1+(1-\lambda)y_2))\\
  	&=&\varphi(\lambda y_1+(1-\lambda)y_2).
  \end{eqnarray*}
  That is,
  \begin{eqnarray*}
  	\varphi(\lambda y_1+(1-\lambda)y_2)<\lambda \varphi(y_1)+(1-\lambda) \varphi (y_2).
  \end{eqnarray*}
  Therefore, according to (ii) of Lemma 3.2, it follows from (\ref{3.24}) and the strict convexity of $\varphi$ that $\varphi$ is $C^1$ on $\Psi_e$ for any $e\in \mathbb{S}^{n-1}$. This implies that $h_K$ is $C^1$ on $\mathbb{R}^{n}\setminus\{o\}$, and $\partial K$ contains no segment. This completes the proof of (i) in Theorem 3.1.
  
  We next prove (ii) in Theorem 3.1. From the assumption that $f$ is continuous, and noting that $\varphi$ is $C^1$ on $\mbox{cl} \Psi_e$ for any $e\in \mathbb{S}^{n-1}$, it follows that the right hand side of (\ref{3.9}) is continuous. Using (i) of Lemma 3.3 and the strict convexity of $\varphi$ on $\Psi_e$, we have that there exists an open ball $B\subset \Psi_e$ centred at $o$ such that $\varphi$ is $C^{1, \alpha}$ on $B$ for any $\alpha\in (0, 1)$. This implies that $h_K$ is locally $C^{1, \alpha}$ on $\mathbb{S}^{n-1}$. Thus it follows from the compactness of $\mathbb{S}^{n-1}$ that $h_K$ is globally $C^{1, \alpha}$ on $\mathbb{S}^{n-1}$, which finishes the proof of (ii) in Theorem 3.1.
  
  Now let us show (iii) of Theorem 3.1. Since $\varphi$ is $C^{1, \alpha}$ on $B$, together with the assumption that $f$ is $C^\alpha$ on $\mathbb{S}^{n-1}$ we obtain that the right hand side of (\ref{3.9}) is $C^\alpha$. According to (ii) of Lemma 3.3, we deduce that $\varphi$ is $C^{2, \alpha}$ on an open ball $\widetilde{B}\subset B$ of $e^\perp$ centred at $o$. This implies that $h_K$ is locally $C^{2, \alpha}$ on $\mathbb{S}^{n-1}$. Thus $h_K$ is globally $C^{2, \alpha}$ on $\mathbb{S}^{n-1}$ based on the compactness of $\mathbb{S}^{n-1}$. This gives the proof of (iii) in Theorem 3.1.\hfill $\square$
  
\section{\bf Existence of symmetric solutions}
\indent In this section, we are going to prove the existence of symmetric solutions to the $L_p$ Gaussian Minkowski problem for $p\leq 0$. The proof is divided into two parts of case $p<0$ and case $p=0$.
\subsection{The proof of case $p<0$}\hfill \break
\indent We shall use a variational method and an approximation argument to solve the $L_p$ Gaussian Minkowski problem for $p< 0$. The crucial step of the variational method is to convert the problem for $p<0$ into an optimization problem whose optimizer is exactly the solution to the original problem. To find such an optimization problem, it is critical that there is certain variational formula that would lead to the probabilistic measure being studied. It is important to note that the variational formula is due to Liu \cite{73} (see Lemma 2.2). At the end of this subsection, an approximation argument will be used to obtain the desired solution.

Let $\mu$ be a non-zero even finite Borel measure on $\mathbb{S}^{n-1}$. For $p<0$ and $f\in C_e^+(\mathbb{S}^{n-1})$, define the functional $\mathcal{J}: C_e^+(\mathbb{S}^{n-1})\rightarrow \mathbb{R}$ by
\begin{eqnarray}\label{4.1z}
\mathcal{J}(f)=-\frac{1}{p}\int_{\mathbb{S}^{n-1}}f(v)^pd\mu(v),
\end{eqnarray}
and define
\begin{eqnarray}\label{4.2z}
	\mathcal{J}(Q)=\mathcal{J}(h_Q)
\end{eqnarray}
for $Q\in \mathcal{K}_e^n$. Our variational method involves the following maximization problem:
\begin{eqnarray*}
\sup\left\lbrace\mathcal{J}(Q): \gamma_n(Q)=\frac{1}{2}~\mbox{and}~ Q\in \mathcal{K}_e^n\right\rbrace.	
\end{eqnarray*}

 {\it \noindent{\bf Lemma 4.1.}~For $p<0$ and $K\in \mathcal{K}_e^n$, let $K$ be a maximizer to the maximization problem
 \begin{eqnarray}\label{4.3z}
 	\sup\left\lbrace\mathcal{J}(Q): \gamma_n(Q)=\frac{1}{2}~\mbox{and}~ Q\in \mathcal{K}_e^n\right\rbrace	
 \end{eqnarray}
if and only if $h_K$ is a maximizer to the maximization problem 
 \begin{eqnarray}\label{4.4z}
 	\sup\left\lbrace\mathcal{J}(f): \gamma_n([f])=\frac{1}{2}~\mbox{and}~ f\in  C_e^+(\mathbb{S}^{n-1})\right\rbrace.	
 \end{eqnarray}}

{\bf Proof.}~For $f\in  C_e^+(\mathbb{S}^{n-1}),$ the Wulff shape $[f]$ generated by $f$ is 
\begin{eqnarray*}
	\left[ f\right]=\bigcap_{v\in \mathbb{S}^{n-1}}\left\lbrace x\in \mathbb{R}^{n}: \langle v, x\rangle\leq f(v)\right\rbrace. 
\end{eqnarray*}
It is easy to see that
\begin{eqnarray}\label{4.5z}
h_{[f]}\leq f~~~~\mbox{and}~~~~[h_{[f]}]=[f].
\end{eqnarray}
Noting $p<0$, we have
\begin{eqnarray}\label{4.6z}
	\mathcal{J}(f)&=&-\frac{1}{p}\int_{\mathbb{S}^{n-1}}f(v)^pd\mu(v) \nonumber \\
	&\leq& -\frac{1}{p}\int_{\mathbb{S}^{n-1}}h_{[f]}(v)^pd\mu(v) \nonumber \\
	&=&\mathcal{J}(h_{[f]}).
\end{eqnarray}
From the second equation of \eqref{4.5z}, it follows that
\begin{eqnarray*}
	\gamma_n([f])=\gamma_n([h_{[f]}]).
\end{eqnarray*}
Since $[f]\in \mathcal{K}_e^n$ for $f\in  C_e^+(\mathbb{S}^{n-1}),$ we get $h_{[f]}\in C_e^+(\mathbb{S}^{n-1})$. Thus,
\begin{eqnarray*}
\mathcal{J}(h_{[f]})\in \left\lbrace\mathcal{J}(f): \gamma_n([f])=\frac{1}{2}~\mbox{and}~ f\in  C_e^+(\mathbb{S}^{n-1})\right\rbrace.	
\end{eqnarray*}
Based on inequality \eqref{4.6z}, we may restrict our attention to the set of all support functions of origin-symmetric convex bodies. Namely,
 \begin{eqnarray}\label{4.7z}
	&&\sup\left\lbrace\mathcal{J}(f): \gamma_n([f])=\frac{1}{2}~\mbox{and}~ f\in  C_e^+(\mathbb{S}^{n-1})\right\rbrace \nonumber\\
	&=&	\sup\left\lbrace\mathcal{J}(h_Q): \gamma_n([h_Q])=\frac{1}{2}~\mbox{and}~ Q\in  \mathcal{K}_e^n\right\rbrace \nonumber \\
	&=&	\sup\left\lbrace\mathcal{J}(Q): \gamma_n(Q)=\frac{1}{2}~\mbox{and}~ Q\in  \mathcal{K}_e^n\right\rbrace.
\end{eqnarray}

On one hand, let $K\in\mathcal{K}_e^n$ be a maximizer to the maximization problem \eqref{4.3z}. Then from \eqref{4.7z} we see 
\begin{eqnarray*}
\sup\left\lbrace\mathcal{J}(f): \gamma_n([f])=\frac{1}{2}~\mbox{and}~ f\in  C_e^+(\mathbb{S}^{n-1})\right\rbrace=\mathcal{J}(K)=\mathcal{J}(h_K).
\end{eqnarray*}
That is, $h_K$ is a maximizer to the maximization problem \eqref{4.4z}. 

On the other hand, let $h_K$ be a maximizer to the maximization problem \eqref{4.4z}
for $K\in  \mathcal{K}_e^n$. Again from \eqref{4.7z}, we have
\begin{eqnarray*}
	\sup\left\lbrace\mathcal{J}(Q): \gamma_n(Q)=\frac{1}{2}~\mbox{and}~ Q\in  \mathcal{K}_e^n\right\rbrace=\mathcal{J}(h_K)=\mathcal{J}(K),
\end{eqnarray*}
i.e., $K$ is a maximizer to the maximization problem \eqref{4.3z}. \hfill $\square$

The following result is to show the existence of a maximizer to the maximization problem \eqref{4.3z} for the case where $p<0$ and the measure $\mu$ is a non-zero even finite Borel measure that vanishes on great subspheres of $\mathbb{S}^{n-1}$.

{\it \noindent{\bf Lemma 4.2.}~Let $p<0$. If $\mu$ is a non-zero even finite Borel measure on $\mathbb{S}^{n-1}$ that vanishes on great subspheres, then there exists a convex body $K\in \mathcal{K}_e^n$ with $\gamma_n(K)=\frac{1}{2}$ such that
\begin{eqnarray*}
	\mathcal{J}(K)=\sup\left\lbrace\mathcal{J}(Q): \gamma_n(Q)=\frac{1}{2}~\mbox{and}~ Q\in  \mathcal{K}_e^n\right\rbrace.
\end{eqnarray*}	
}

{\bf Proof.}~Let $\left\lbrace Q_i \right\rbrace \subset \mathcal{K}_e^n$ be a maximizing sequence, i.e.,
\begin{eqnarray*}
\lim_{i\rightarrow \infty}\mathcal{J}(Q_i)=\sup\left\lbrace\mathcal{J}(Q): \gamma_n(Q)=\frac{1}{2}~\mbox{and}~ Q\in  \mathcal{K}_e^n\right\rbrace.
\end{eqnarray*}	
Choosing a constant $c_0>0$ such that $\gamma_n(c_0 B)=\frac{1}{2}$ for the unit ball $B$ in $\mathbb{R}^n$, it is obvious that 
\begin{eqnarray}\label{4.8z}
\mathcal{J}(c_0 B)\in\left\lbrace\mathcal{J}(Q): \gamma_n(Q)=\frac{1}{2}~\mbox{and}~ Q\in  \mathcal{K}_e^n\right\rbrace.
\end{eqnarray}	
Thus, it follows that for sufficiently large $i$,
\begin{eqnarray}\label{4.9z}
	\mathcal{J}(Q_i)\geq\mathcal{J}(c_0 B)=-\frac{1}{p}\int_{\mathbb{S}^{n-1}}c_0^pd\mu(v)=-\frac{1}{p}c_0^p\mu(\mathbb{S}^{n-1})>0.
\end{eqnarray}	

We first show that $h_{Q_i}$ has the uniformly positive bound from below, proceeding by contradiction. Suppose that there exists a subsequence, denoted also by $Q_i$ and $\gamma_n(Q_i)=\frac{1}{2}$, which converges to a compact convex subset $Q_0\subset \mathbb{R}^n$ with $h_{Q_0}(v_0)=0$ for some $v_0\in \mathbb{S}^{n-1}$. Let $h_{Q_i}(v_i)=\min_{v\in\mathbb{S}^{n-1}}h_{Q_i}(v)$ for $v_i\in \mathbb{S}^{n-1}$. Then 
\begin{eqnarray*}
	\lim_{i\rightarrow \infty}h_{Q_i}(v_i)=h_{Q_0}(v_0)=0.
\end{eqnarray*}	
Since 
\begin{eqnarray*}
	Q_i\subset \left\lbrace x\in \mathbb{R}^n: |\langle x, v_i\rangle|\leq h_{Q_i}(v_i)\right\rbrace, 
\end{eqnarray*}	
we have
\begin{eqnarray*}
	\gamma_n(Q_i)\leq\gamma_n\left(  \left\lbrace x\in \mathbb{R}^n: |\langle x, v_i\rangle|\leq h_{Q_i}(v_i)\right\rbrace\right)\rightarrow 0, 
\end{eqnarray*}	
which contradicts to $\gamma_n(Q_i)=\frac{1}{2}$.

Next, we will show that $Q_i$ is uniformly bounded, also arguing by contradiction. If not, there exists a subsequence, denoted again by $Q_i$, such that $Q_i^*$ converges to $Q_0^*$ with $h_{Q^*_0}(\pm v_0)=0$ for some $v_0\in \mathbb{S}^{n-1}$. For each $\delta >0$, define
\begin{eqnarray*}
	\Omega_\delta(v_0)=\left\lbrace v\in \mathbb{S}^{n-1}: |\langle v, v_0\rangle|>\delta\right\rbrace. 
\end{eqnarray*}	
Since $h_{Q^*_0}(v_0)=0$ and $Q_i^*\rightarrow Q_0^*$, we have $\lim_{i\rightarrow \infty}h_{Q_i^*}(v_0)=0$. Thus, $\rho_{Q_i^*}\rightarrow 0$ uniformly on $\Omega_\delta(v_0)$. By hypothesis, $\mu$ vanishes on all great subspheres. Then
\begin{eqnarray*}
\mu\left( \mathbb{S}^{n-1}\cap v_0^\bot\right)=0,
\end{eqnarray*}	
where $v_0^\bot$ is the co-dimension $1$ subspace. Hence,
\begin{eqnarray*}
\lim_{\delta\rightarrow 0^+}\mu\left( \mathbb{S}^{n-1}\backslash \Omega_\delta(v_0)\right)=\mu\left( \mathbb{S}^{n-1}\cap v_0^\bot\right)=0.
\end{eqnarray*}	
Since $h_{Q_i}$ has the uniformly positive bound from below, there exists a constant $R>0$ such that $Q_i^*\subset RB$ for the unit ball $B$ in $\mathbb{R}^n$. For $p<0$, we have
\begin{eqnarray}\label{4.10z}
\int_{\mathbb{S}^{n-1}}h_{Q_i}^p(v)d\mu(v)&=&\int_{\Omega_\delta(v_0)}\rho_{Q^*_i}^{-p}(v)d\mu(v)+\int_{\mathbb{S}^{n-1}\setminus\Omega_\delta(v_0)}\rho_{Q^*_i}^{-p}(v)d\mu(v) \nonumber \\
&\leq&\int_{\Omega_\delta(v_0)}\rho_{Q^*_i}^{-p}(v)d\mu(v)+R^{-p}\mu\left(\mathbb{S}^{n-1}\setminus\Omega_\delta(v_0) \right).
\end{eqnarray}	
For each $\varepsilon>0$, we choose $\delta_0>0$ such that
\begin{eqnarray}\label{4.11z}
R^{-p}\mu\left(\mathbb{S}^{n-1}\setminus\Omega_{\delta_0}(v_0) \right)<\frac{\varepsilon}{2}.
\end{eqnarray}	
Since $\rho_{Q_i^*}\rightarrow 0$ uniformly on $\Omega_{\delta_0}(v_0)$ and $p<0$, we can find an $i_0$ such that for all $i\geq i_0$ and each $v\in \Omega_{\delta_0}(v_0)$,
\begin{eqnarray}\label{4.12z}
\mu(\mathbb{S}^{n-1})\rho_{Q_i^*}(v)^{-p}<\frac{\varepsilon}{2}.
\end{eqnarray}	
Together \eqref{4.10z}, \eqref{4.11z} and \eqref{4.12z}, it follows that
\begin{eqnarray*}
\int_{\mathbb{S}^{n-1}}h_{Q_i}^p(v)d\mu(v) \leq \varepsilon.
\end{eqnarray*}	
Namely, 
\begin{eqnarray*}
	\lim_{i\rightarrow \infty}\int_{\mathbb{S}^{n-1}}h_{Q_i}^p(v)d\mu(v)=0.
\end{eqnarray*}	
From the definition \eqref{4.1z} of $\mathcal{J}$ and \eqref{4.9z}, we obtain a contradiction. Therefore, the sequence $h_{Q_i}$
 has uniformly positive bounds from below and above.

By Blaschke's selection theorem, we obtain $Q_i$ subsequently converges to a convex body $K\in \mathcal{K}^n_e$ with $\gamma_n(K)=\frac{1}{2}$ such that
\begin{eqnarray*}
	\mathcal{J}(K)=\sup\left\lbrace\mathcal{J}(Q): \gamma_n(Q)=\frac{1}{2}~\mbox{and}~ Q\in  \mathcal{K}_e^n\right\rbrace.
\end{eqnarray*}	
\hfill $\square$

When the even measure $\mu$ vanishes on all great subspheres of $\mathbb{S}^{n-1}$, the following theorem provides the existence of a solution to  the $L_p$ Gaussian Minkowski problem for $p<0$.

{\it \noindent{\bf Theorem 4.1.}~Let $p< 0$ and $\mu$ be a non-zero even finite Borel measure on $\mathbb{S}^{n-1}$. If $\mu$ vanishes on all great subspheres, then there exists an origin-symmetric convex body $K\in \mathcal{K}_e^n$ with $\gamma_n(K)=\frac{1}{2}$ such that
	\begin{eqnarray*}
		S_{p, \gamma_n}(K, \cdot)=\frac{S_{p, \gamma_n}(K, \mathbb{S}^{n-1})}{\mu(\mathbb{S}^{n-1})}\mu.
\end{eqnarray*}} 

{\bf Proof.}~By Lemma 4.2, there exists a convex body $K\in \mathcal{K}_e^n$ with $\gamma_n(K)=\frac{1}{2}$ such that
\begin{eqnarray*}
	\mathcal{J}(K)=\sup\left\lbrace\mathcal{J}(Q): \gamma_n(Q)=\frac{1}{2}~\mbox{and}~ Q\in  \mathcal{K}_e^n\right\rbrace.
\end{eqnarray*}	
According to Lemma 4.1, this implies
\begin{eqnarray*}
	\mathcal{J}(h_K)=\sup\left\lbrace\mathcal{J}(f): \gamma_n([f])=\frac{1}{2}~\mbox{and}~ f\in  C_e^+(\mathbb{S}^{n-1})\right\rbrace.	
\end{eqnarray*}
For any $g\in C(\mathbb{S}^{n-1})$, $\lambda \in \mathbb{R}^n$, and $t\in (-\delta, \delta)$ where $\delta$ sufficiently small, define
\begin{eqnarray*}
	\mathscr{L}_1(t, \lambda)=\mathcal{J}(h_t)+\lambda \left( \gamma_n([h_t])-\frac{1}{2}\right),
\end{eqnarray*}
where $h_t$ is given by
\begin{eqnarray*}
h_t(v)=(h_K(v)^p+tg(v)^p)^{\frac{1}{p}}
\end{eqnarray*}
for $v\in \mathbb{S}^{n-1}$ and $p<0$. Hence, it follows from the method of Lagrange multipliers that
\begin{eqnarray*}
\left.\frac{\partial}{\partial t}\mathscr{L}_1(t, \lambda) \right|_{t=0}=0.
\end{eqnarray*}
From the definition \eqref{4.1z} of $\mathcal{J}$ and Lemma 2.2, we see
\begin{eqnarray*}
	0&=&\left.\frac{\partial}{\partial t}\left( \mathcal{J}(h_t)+\lambda \left( \gamma_n([h_t])-\frac{1}{2}\right)\right)  \right|_{t=0}\\
	&=&-\frac{1}{p}\int_{\mathbb{S}^{n-1}}g(v)^pd\mu(v)+\frac{\lambda}{p}\int_{\mathbb{S}^{n-1}}g(v)^pdS_{p, \gamma_n}(K, v).
\end{eqnarray*}
Namely,
\begin{eqnarray}\label{4.13z}
\int_{\mathbb{S}^{n-1}}g(v)^pd\mu(v)=\lambda\int_{\mathbb{S}^{n-1}}g(v)^pdS_{p, \gamma_n}(K, v).
\end{eqnarray}
Thus, by the arbitrariness of $g$, we have
\begin{eqnarray}\label{4.14z}
\mu=\lambda S_{p, \gamma_n}(K, \cdot).
\end{eqnarray}
Letting $g\equiv 1$ in \eqref{4.13z}, we have
\begin{eqnarray*}
\mu(\mathbb{S}^{n-1})=\int_{\mathbb{S}^{n-1}}d\mu(v)=\lambda\int_{\mathbb{S}^{n-1}}d S_{p, \gamma_n}(K, v)=\lambda S_{p, \gamma_n}(K, \mathbb{S}^{n-1}).
\end{eqnarray*}
That is,
\begin{eqnarray}\label{4.15z}
	\lambda=\frac{\mu(\mathbb{S}^{n-1})}{ S_{p, \gamma_n}(K, \mathbb{S}^{n-1})}
\end{eqnarray}
Combining \eqref{4.14z} with \eqref{4.15z}, we obtain
	\begin{eqnarray*}
	S_{p, \gamma_n}(K, \cdot)=\frac{S_{p, \gamma_n}(K, \mathbb{S}^{n-1})}{\mu(\mathbb{S}^{n-1})}\mu.
\end{eqnarray*}
\hfill $\square$

In the following we solve the even $L_p$ Gaussian Minkowski problem for $p<0$ in its full generality, which is obtained by using Theorem 4.1 and an approximation argument. Namely, we will complete the proof of case $p<0$ in Theorem 1.1.

{\bf Proof of Theorem 1.1 ($p<0$).}~~Let $f_j$ be a sequence of positive and smooth even functions on $\mathbb{S}^{n-1}$ such that $\{\mu_j\}$, $d\mu_j=f_j(v)dv$, converges to $\mu$ weakly. Then, it is easy to see that $\mu_j$ satisfies the conditions in Theorem 4.1. In particular, $\mu_j$ vanishes on all great subspheres of $\mathbb{S}^{n-1}$. Hence, it follows from Theorem 4.1 that there exists a weak solution $K_j\in \mathcal{K}_e^n$ with $\gamma_n(K_j)=\frac{1}{2}$ such that
	\begin{eqnarray}\label{4.16z}
	S_{p, \gamma_n}(K_j, \cdot)=\frac{S_{p, \gamma_n}(K_j, \mathbb{S}^{n-1})}{\mu_j(\mathbb{S}^{n-1})}\mu_j.
\end{eqnarray}
Thus, $h_{K_j}$ satisfies the following Monge-Amp\`{e}re type equation on $\mathbb{S}^{n-1}$ in sense of measure:
 \begin{eqnarray}\label{4.17z}
	\frac{1}{(\sqrt{2\pi})^n}h^{1-p}_{K_j}e^{-\frac{\left|\nabla h_{K_j}\right|^2+h^2_{K_j} }{2}}\det\left(\nabla^2h_{K_j}+h_{K_j}I\right)=\frac{S_{p, \gamma_n}(K_j, \mathbb{S}^{n-1})}{\mu_j(\mathbb{S}^{n-1})}f_j.
\end{eqnarray} 
By (iii) of Theorem 3.1, we see that $h_{K_j}$ belongs to $C^{2, \alpha}(\mathbb{S}^{n-1})$, provided $0<f_j\in C^{\alpha}(\mathbb{S}^{n-1})$ for $\alpha\in (0, 1)$. Thus, according to the standard regularity theory of Monge-Amp\`{e}re equation, we know that $K_j$ is smooth and uniformly convex when $f_j$ is smooth and positive on $\mathbb{S}^{n-1}$.

Let $h_{K_j}(v_j)=\max_{v\in\mathbb{S}^{n-1}}h_{K_j}(v)$ for some $v_j\in \mathbb{S}^{n-1}$. Since $\mu$ is not concentrated on any great subspheres of $\mathbb{S}^{n-1}$, there exist $\delta, \zeta>0$ such that for any $u\in \mathbb{S}^{n-1}$,
\begin{eqnarray*}
\int_{\{v\in\mathbb{S}^{n-1}:|\langle v, u\rangle|>\zeta\}}d\mu(v)\geq \delta>0.
\end{eqnarray*}
Hence, for enough large $j$,
\begin{eqnarray}\label{4.18z}
	\int_{\{v\in\mathbb{S}^{n-1}:|\langle v, v_j\rangle|>\zeta\}}f_j(v)dv\geq \delta>0.
\end{eqnarray}
Noting $\rho_{K_j}(v_j)=h_{K_j}(v_j)$, we get that for any $v\in \mathbb{S}^{n-1}$, 
\begin{eqnarray}\label{4.19z}
h_{K_j}(v)\geq |\langle v, v_j\rangle|\rho_{K_j}(v_j)=|\langle v, v_j\rangle|h_{K_j}(v_j).
\end{eqnarray}
Note also that
 \begin{eqnarray}\label{4.20z}
	\det\left(\nabla^2h_{K_j}+h_{K_j}I\right)=\frac{\left( h_{K_j}^2+|\nabla h_{K_j}|^2\right)^{\frac{n}{2}}}{h_{K_j}},
\end{eqnarray} 
and for $u, v\in \mathbb{S}^{n-1}$,
 \begin{eqnarray}\label{4.21z}
\rho_{K_j}^2(u)=h_{K_j}^2(v)+|\nabla h_{K_j}(v)|^2,
\end{eqnarray} 
where $u$ and $v$ are associated by
 \begin{eqnarray*}
	\rho_{K_j}(u)u=h_{K_j}(v)v+\nabla h_{K_j}(v).
\end{eqnarray*} 
Thus, it follows from \eqref{4.17z}, \eqref{4.18z}, \eqref{4.19z}, \eqref{4.20z} and \eqref{4.21z} that
\begin{eqnarray*}
0&<& \frac{S_{p, \gamma_n}(K_j, \mathbb{S}^{n-1})}{\mu_j(\mathbb{S}^{n-1})}\delta\\
&\leq& \frac{S_{p, \gamma_n}(K_j, \mathbb{S}^{n-1})}{\mu_j(\mathbb{S}^{n-1})}\int_{\{v\in\mathbb{S}^{n-1}:|\langle v, v_j\rangle|>\zeta\}}f_j(v)dv
\end{eqnarray*}
\begin{eqnarray*}
&=& \frac{1}{(\sqrt{2\pi})^n}\int_{\{v\in\mathbb{S}^{n-1}:|\langle v, v_j\rangle|>\zeta\}}e^{-\frac{\left|\nabla h_{K_j}(v)\right|^2+h^2_{K_j} (v)}{2}}h^{1-p}_{K_j}(v)\det\left(\nabla^2h_{K_j}(v)+h_{K_j}(v)I\right)dv\\
&\leq &\frac{1}{(\sqrt{2\pi})^n}e^{-\frac{\left( \zeta h_{K_j} (v_j)\right)^2}{2}}h^{1-p}_{K_j}(v_j)\int_{\{v\in\mathbb{S}^{n-1}:|\langle v, v_j\rangle|>\zeta\}}\det\left(\nabla^2h_{K_j}(v)+h_{K_j}(v)I\right)dv\\
&= &\frac{1}{(\sqrt{2\pi})^n}\frac{h^{1-p}_{K_j}(v_j)}{e^{\frac{\left( \zeta h_{K_j} (v_j)\right)^2}{2}}}\int_{\{v\in\mathbb{S}^{n-1}:|\langle v, v_j\rangle|>\zeta\}}\frac{\left( h_{K_j}^2(v)+|\nabla h_{K_j}(v)|^2\right)^{\frac{n}{2}}}{h_{K_j}(v)}dv\\
&\leq &\frac{1}{(\sqrt{2\pi})^n}\frac{h^{-p}_{K_j}(v_j)}{\zeta e^{\frac{\left( \zeta h_{K_j} (v_j)\right)^2}{2}}}\int_{\{v\in\mathbb{S}^{n-1}:|\langle v, v_j\rangle|>\zeta\}}\left( h_{K_j}^2(v)+|\nabla h_{K_j}(v)|^2\right)^{\frac{n}{2}}dv\\
&\leq &\frac{1}{(\sqrt{2\pi})^n}\frac{h^{-p}_{K_j}(v_j)}{\zeta e^{\frac{\left( \zeta h_{K_j} (v_j)\right)^2}{2}}}\int_{\mathbb{S}^{n-1}}\rho_{K_j}^n(u)du\\
&\leq &\frac{1}{(\sqrt{2\pi})^n}\frac{h^{n-p}_{K_j}(v_j)}{\zeta e^{\frac{\left( \zeta h_{K_j} (v_j)\right)^2}{2}}}\mathcal{H}^{n-1}(\mathbb{S}^{n-1})\rightarrow 0,
\end{eqnarray*}
as $h_{K_j}(v_j)\rightarrow \infty$. This is a contradiction. Therefore, $K_j$ is uniformly bounded. The uniformly positive lower bound for $h_{K_j}$ is garanteed by $\gamma_n(K_j)=\frac{1}{2}$. Its proof has been given in Lemma 4.2. Thus, there exists a positive constant $C_1$ independent of $j$ such that 
\begin{eqnarray*}
	\frac{1}{C_1}\leq h_{K_j}\leq C_1.
\end{eqnarray*} 

By Blaschke's selection theorem, we know that there exists a subsequence $K_i$, of $K_j$, with $\gamma_n(K_i)=\frac{1}{2}$ such that $K_i\rightarrow K\in \mathcal{K}^n_e$. Thus, it follows from the weak convergence of the $L_p$ Gaussian surface area measure that \eqref{4.16z} implies
\begin{eqnarray*}
	S_{p, \gamma_n}(K, \cdot)=\frac{S_{p, \gamma_n}(K, \mathbb{S}^{n-1})}{\mu(\mathbb{S}^{n-1})}\mu.
\end{eqnarray*}
and $\gamma_n(K)=\frac{1}{2}.$
\hfill $\square$
\subsection{The proof of case $p=0$}\hfill \break

In this subsection, we show the existence of symmetric solutions to the $L_0$ Gaussian Minkowski problem by studying the related Monge-Amp\`{e}re type functional on the unit sphere.

Let $K\in \mathcal{K}_o^n$. Given any $f\in C(\mathbb{S}^{n-1})$, there is a $\delta>0$ such that $h_K(v)+tf(v)>0$ for all $v\in \mathbb{S}^{n-1}$ and $t\in (-\delta, \delta)$. Defining $h_t(v)=h_K(v)+tf(v)$, we consider a family of convex bodies
\begin{eqnarray}\label{4.22z}
	\left[ h_t\right]=\bigcap_{v\in \mathbb{S}^{n-1}}\left\lbrace x\in \mathbb{R}^{n}: \langle v, x\rangle\leq h_t(v)\right\rbrace. 
\end{eqnarray}
Let $h(\cdot, t)$ and $\rho(\cdot, t)$ be respectively the support function and radial function, on $\mathbb{S}^{n-1}$, of $[h_t]$. That is, $h(\cdot, t)=h_{[h_t]}$ and $\rho(\cdot, t)=\rho_{[h_t]}$. Then the following lemma was established in \cite{1b}.
 
 {\it \noindent{\bf Lemma 4.3 \cite{1b}.}~Suppose that $\partial K$ is $C^1$ and strictly convex at $x\in \partial K$ for $K\in \mathcal{K}_o^n$. Then the limits
 	\begin{eqnarray*}
h'(v)=\lim_{t\rightarrow 0}\frac{h(v,t)-h(v, 0)}{t}, ~~~\rho'(u) =\lim_{t\rightarrow 0}\frac{\rho(u,t)-\rho(u, 0)}{t}
\end{eqnarray*}	
exist, where $v$ is the unit outer normal of $\partial K$ at $x$ and $u=\frac{x}{|x|}=\alpha_K^*(v)$. Moreover,	
	\begin{eqnarray}\label{4.23z}
 h'(v) =f(v),
\end{eqnarray}	
	\begin{eqnarray}\label{4.24z}
	\frac{\rho'(u)}{\rho_K(u)} =\frac{h'(v)}{h_K(v)}.
\end{eqnarray} }

{\bf Proof of Theorem 1.1 ($p=0$).}~~We will divide into four steps to complete the proof.

{\bf Step 1.}~~We first show the following variational formula for the Gaussian volume. Namely, if $\partial K$ is $C^1$ and strictly convex for $K\in \mathcal{K}_o^n$, we prove	
\begin{eqnarray}\label{4.25z}
\left. \frac{d \gamma_n([h_t])}{dt}\right| _{t=0}=\frac{1}{\left(\sqrt{2\pi} \right)^n}\int_{\mathbb{S}^{n-1}}e^{-\frac{\rho_K(u)^2}{2}}\rho_K(u)^n\frac{f(\alpha_K(u))}{h_K(\alpha_K(u))}du.
\end{eqnarray}

Using the polar coordinates formula, it follows that
\begin{eqnarray*}
\gamma_n([h_t])=\frac{1}{\left(\sqrt{2\pi} \right)^n}\int_{\mathbb{S}^{n-1}}\int_0^{\rho_{[h_t]}(u)}e^{-\frac{r^2}{2}}r^{n-1}drdu.
\end{eqnarray*}
Noting $K\in \mathcal{K}_o^n$ and $h_t\rightarrow h_K$ uniformly on $\mathbb{S}^{n-1}$ as $t\rightarrow 0$, it follows from the Aleksandrov convergence lemma \cite{77} that
\begin{eqnarray*}
\lim_{t\rightarrow 0}[h_t]=[h_K]=K.
\end{eqnarray*}
From this, we obtain that there exists $M_1>0$ such that $[h_t]\subset M_1 B$ for the unit ball $B$ in $\mathbb{R}^n$. Let
\begin{eqnarray*}
	\Gamma (s)=\int_0^se^{-\frac{r^2}{2}}r^{n-1}dr.
\end{eqnarray*}
Then by mean value theorem
\begin{eqnarray}\label{4.26z}
\left| \Gamma(\rho_{[h_t]}(u))-\Gamma(\rho_K(u)) \right|=\left| \Gamma'(\xi)\right|  \left|\rho_{[h_t]}(u)-\rho_K(u)\right|, 
\end{eqnarray}
where $\xi \in (0, M_1]$. From the definition of $\Gamma$, there exists $M_2>0$ such that 
\begin{eqnarray}\label{4.27z}
	\left| \Gamma'(\xi) \right| \leq M_2.
\end{eqnarray}
By \eqref{4.24z}, we have that for all $u\in \mathbb{S}^{n-1}$,
\begin{eqnarray*}
\lim_{t\rightarrow 0}\frac{\rho_{[h_t]}(u)-\rho_K(u)}{t}=\frac{f(\alpha_K(u))}{h_K(\alpha_K(u))}\rho_K(u).
\end{eqnarray*}
Thus, it follows that there exists $M_3>0$ such that when $t\rightarrow 0$,
\begin{eqnarray}\label{4.28z}
\left|\rho_{[h_t]}(u)-\rho_K(u) \right|\leq M_3 |t|,
\end{eqnarray}
for all $u\in \mathbb{S}^{n-1}$. Therefore, from \eqref{4.26z}, \eqref{4.27z} and \eqref{4.28z} we obtain that
for all $u\in \mathbb{S}^{n-1}$,
\begin{eqnarray*}
	\left| \Gamma(\rho_{[h_t]}(u))-\Gamma(\rho_K(u)) \right|\leq M_2M_3|t|. 
\end{eqnarray*}
Hence, it follows from Lemma 4.3 and the dominated convergence theorem that
\begin{eqnarray*}
	\lim_{t\rightarrow 0}\frac{\gamma_n([h_t])-\gamma_n(K)}{t}=\frac{1}{\left(\sqrt{2\pi} \right)^n}\int_{\mathbb{S}^{n-1}}e^{-\frac{\rho_K(u)^2}{2}}\rho_K(u)^n\frac{f(\alpha_K(u))}{h_K(\alpha_K(u))}du.
\end{eqnarray*}
This completes Step 1.

Let $d\mu =g(v)dv$ for $v\in \mathbb{S}^{n-1}$, where $g$ is an even function on $\mathbb{S}^{n-1}$ with $1/C\leq g\leq C$
for some constant $C>0$. For $K\in \mathcal{K}_o^n$, define
\begin{eqnarray}\label{4.29z}
	\mathcal{E}(K)=-\int_{\mathbb{S}^{n-1}}\log h_K(v)d\mu(v).
\end{eqnarray}
Then we consider the following maximization problem:
\begin{eqnarray}\label{4.30z}
\sup\left\lbrace \mathcal{E}(Q):\gamma_n(Q)=\frac{1}{2}~~\mbox{and}~~Q\in \mathcal{K}_e^n\right\rbrace.
\end{eqnarray}

{\bf Step 2.}~~We show the existence of an optimizer for the maximization problem \eqref{4.30z}. Namely, there exists a convex body $K\in \mathcal{K}_e^n$ with $\gamma_n(K)=\frac{1}{2}$ such that
\begin{eqnarray}\label{4.31z}
	\mathcal{E}(K)=\sup\left\lbrace \mathcal{E}(Q):\gamma_n(Q)=\frac{1}{2}~~\mbox{and}~~Q\in \mathcal{K}_e^n\right\rbrace.
\end{eqnarray}
Moreover, $\partial K$ is $C^1$ and strictly convex.

Let
\begin{eqnarray*}
	\lim_{i\rightarrow \infty}\mathcal{E}(Q_i)=\sup\left\lbrace\mathcal{E}(Q): \gamma_n(Q)=\frac{1}{2}~\mbox{and}~ Q\in  \mathcal{K}_e^n\right\rbrace,
\end{eqnarray*}	
where $Q_i\in \mathcal{K}^n_e$ and $\gamma_n(Q_i)=\frac{1}{2}$. We claim that $Q_i$
is uniformly bounded and argue by contradiction. Assume
\begin{eqnarray*}
R_i=\max_{u\in\mathbb{S}^{n-1}}\rho_{Q_i}(u)=\rho_{Q_i}(u_i),
\end{eqnarray*}	
where $u_i$ is one of the unit vectors at which the maximum occurs. Since $Q_i\in \mathcal{K}^n_e$ and $\left\lbrace  ru_i:-R_i\leq r\leq R_i \right\rbrace \subset Q_i $, we have that for all $v\in \mathbb{S}^{n-1}$,
\begin{eqnarray*}
h_{Q_i}(v)\geq R_i|\langle u_i, v\rangle|.
\end{eqnarray*}	
Thus, by definition \eqref{4.29z} we see
\begin{eqnarray*}
	\mathcal{E}(Q_i)&=&-\int_{\mathbb{S}^{n-1}}\log h_{Q_i}(v)d\mu(v)\\
	&\leq& -\int_{\mathbb{S}^{n-1}}\log R_ig(v)dv-\int_{\mathbb{S}^{n-1}}\log |\langle u_i, v\rangle| g(v)dv\\
	&\leq& -\frac{\mathcal{H}^{n-1}(\mathbb{S}^{n-1})}{C}\log R_i-\frac{1}{C}\int_{\mathbb{S}^{n-1}}\log|\langle u_i, v\rangle|dv,
\end{eqnarray*}
where the constant $C$ is from the upper bound of $g$. Note that $\log|\langle u_i, v\rangle|$ is an integral function on $\mathbb{S}^{n-1}$ and its integral on $\mathbb{S}^{n-1}$ is independent of $u_i$. Hence, if $R_i\rightarrow \infty$,
\begin{eqnarray*}
	\lim_{i\rightarrow \infty}\mathcal{E}(Q_i)=-\infty.
\end{eqnarray*}	
This cannot occur as $Q_i$ is a maximising sequence. The positive lower bound of $h_{Q_i}$ is guaranteed by
\begin{eqnarray*}
\gamma_n(Q_i)=\frac{1}{2}.
\end{eqnarray*}	
See the previous subsection for its proof. In a word, there exists a constant $C_1>0$ such that for all $v\in \mathbb{S}^{n-1}$,
\begin{eqnarray*}
\frac{1}{C_1}\leq h_{Q_i}(v)\leq C_1.
\end{eqnarray*}	

By Blaschke's selection theorem, the sequence $Q_i$ has a convergent subsequence the limit of which is, say, $K\in \mathcal{K}^n_e$ satisfying $\gamma_n(K)=\frac{1}{2}$. Moreover,
\begin{eqnarray*}
	\mathcal{E}(K)=\sup\left\lbrace \mathcal{E}(Q):\gamma_n(Q)=\frac{1}{2}~~\mbox{and}~~Q\in \mathcal{K}_e^n\right\rbrace.
\end{eqnarray*}

Further, we also show that $\partial K$ is $C^1$ and strictly convex for the maximiser  $K\in \mathcal{K}^n_e$.

For $f\in C(\mathbb{S}^{n-1})$ and $t\in (-\delta, \delta)$ where $\delta>0$ is sufficiently small, define
\begin{eqnarray}\label{4.32z}
\tau_t(v)=h_K(v)e^{tf(v)},
\end{eqnarray}
for all $v\in \mathbb{S}^{n-1}$. Then the Wulff shape associated with $\tau_t$ is
\begin{eqnarray*}
	\left[ \tau_t\right]=\bigcap_{v\in \mathbb{S}^{n-1}}\left\lbrace x\in \mathbb{R}^{n}: \langle v, x\rangle\leq \tau_t(v)\right\rbrace. 
\end{eqnarray*}
Let $h(\cdot, t)$ and $\rho(\cdot, t)$ be respectively the support function and radial function, on $\mathbb{S}^{n-1}$, of $[\tau_t]$. From \eqref{4.32z}, we have
\begin{eqnarray*}
\log \tau_t(v)=\log h_K(v)+tf(v).
\end{eqnarray*}
Thus, it follows by \eqref{2.8z}, \eqref{2.12z} and  \eqref{2.13z} that
	\begin{eqnarray}\label{4.33z}
	\lim_{t\rightarrow 0}\frac{\rho_{[\tau_t]}(u)-\rho_K(u)}{t}=f(\alpha_K(u))\rho_K(u)
\end{eqnarray} 
holds for almost all $u\in \mathbb{S}^{n-1}$, and there exist $\delta_0>0$ and $M>0$ such that
\begin{eqnarray*}
	\left|\rho_{[\tau_t]}(u)-\rho_K(u) \right|\leq M |t|,
\end{eqnarray*}
for all $u\in \mathbb{S}^{n-1}$ and $t\in(-\delta_0, \delta_0)$. Therefore, using the dominated convergence theorem, it follows from \eqref{2.10z}, \eqref{4.33z} and \eqref{1.1a} that
\begin{eqnarray}\label{4.34z}
	\lim_{t\rightarrow 0}\frac{\gamma_n([\tau_t])-\gamma_n(K)}{t}=\int_{\mathbb{S}^{n-1}}f(v)dS_{0, \gamma_n}(K, v).
\end{eqnarray}

Clearly, $h(\cdot, t)\leq e^{tf}h_K$. From this, we see that for $v\in \mathbb{S}^{n-1}$,
\begin{eqnarray}\label{4.35z}
	\lim_{t\rightarrow 0^+}\frac{h(v, t)-h_K(v)}{t}\leq f(v)h_K(v).
\end{eqnarray}
Since $K$ is the maximizer of \eqref{4.30z}, it follows from \eqref{4.34z} and \eqref{4.35z} that for $\lambda\in \mathbb{R}$,
\begin{eqnarray*}
	0&\geq& \left. \frac{d}{dt}\left(\mathcal{E} \left( [\tau_t]\right)+\lambda\left(\gamma_n([\tau_t])-\frac{1}{2} \right)  \right) \right|_{t=0^+}\\
	&\geq& -\int_{\mathbb{S}^{n-1}}f(v)d\mu(v)+\lambda\int_{\mathbb{S}^{n-1}}f(v)dS_{0, \gamma_n}(K, v).
\end{eqnarray*}
Namely, 
\begin{eqnarray}\label{4.36z}
\int_{\mathbb{S}^{n-1}}f(v)d\mu(v)\geq\lambda\int_{\mathbb{S}^{n-1}}f(v)dS_{0, \gamma_n}(K, v).
\end{eqnarray}

From \eqref{1.1a} and \eqref{2.9z}, we deduce
\begin{eqnarray}\label{4.37z}
\ \ \ \ \ \int_{\mathbb{S}^{n-1}}f(v)dS_{0, \gamma_n}(K, v)&=&\frac{1}{\left( \sqrt{2\pi}\right) ^n}\int_{\partial K} f(\nu_K(x))e^{-\frac{|x|^2}{2}}\langle x, \nu_K(x)\rangle d\mathcal{H}^{n-1}(x) \nonumber \\
&=&\frac{1}{\left( \sqrt{2\pi}\right) ^n}\int_{\mathbb{S}^{n-1}} f(v)e^{-\frac{\rho^2_K(\alpha^*_K(v))}{2}}h_K(v)dS(K, v).
\end{eqnarray}
By the arbitrariness of $f$, this together with inequality \eqref{4.36z} gives that for all $v\in\mathbb{S}^{n-1}$,
\begin{eqnarray}\label{4.38z}
	\frac{dS(K, v)}{dv}\leq \frac{ \left( \sqrt{2\pi}\right) ^ng(v)e^{\frac{\rho^2_K(\alpha^*_K(v))}{2}}}{\lambda h_K(v)}.
\end{eqnarray}
Taking $f\equiv 1$ in \eqref{4.36z}, we also obtain
\begin{eqnarray}\label{4.39z}
	\lambda \leq\frac{\mu(\mathbb{S}^{n-1})}{S_{0, \gamma_n}(K, \mathbb{S}^{n-1})}.
\end{eqnarray}

Let $K^*$ be the polar body of the maximiser $K\in \mathcal{K}^n_e$ in \eqref{4.30z}. Then we know that $\rho_{K^*}=1/h_K$. For $f\in C(\mathbb{S}^{n-1})$, define
\begin{eqnarray*}
\tau_t^*=e^{tf}\rho_{K^*}.
\end{eqnarray*}
The convex hull $\langle \tau_t^*\rangle$ generated by $\tau_t^*$ is
\begin{eqnarray*}
	\langle \tau_t^*\rangle=\mbox{conv}\left\lbrace e^{tf(v)}\rho_{K^*}(v)v: v\in \mathbb{S}^{n-1} \right\rbrace. 
\end{eqnarray*}
According to \eqref{2.8z}, we see 
\begin{eqnarray}\label{4.40z}
	\langle \tau_t^*\rangle^*=[\tau_t].
\end{eqnarray}
Denote by $h^*(\cdot, t)$ and $\rho^*(\cdot, t)$ the support and radial functions of $\langle \tau_t^*\rangle$. Then
\begin{eqnarray*}
h^*(\cdot, t)\geq e^{tf}\rho_{K^*}.
\end{eqnarray*}
Thus, from \eqref{4.40z} we see
\begin{eqnarray*}
	h(\cdot, t)\leq e^{-tf}h_K.
\end{eqnarray*}
Therefore, 
\begin{eqnarray}\label{4.41z}
	\lim_{t\rightarrow 0^+}\frac{h(v, t)-h_K(v)}{t}\leq -f(v)h_K(v)
\end{eqnarray}
for $v\in \mathbb{S}^{n-1}$. Since
\begin{eqnarray*}
	\log \tau^*_t(v)=\log \rho_{K^*}(v)+tf(v),
\end{eqnarray*}
it follows from \eqref{2.14z} that 
	\begin{eqnarray}\label{4.42z}
	\lim_{t\rightarrow 0}\frac{\rho_{[\tau_t]}(u)-\rho_K(u)}{t}=-f(\alpha_K(u))\rho_K(u).
\end{eqnarray} 
Besides, from \eqref{2.15z} there exist $\delta_0>0$ and $M>0$ such that
\begin{eqnarray*}
	\left|\rho_{[\tau_t]}(u)-\rho_K(u) \right|\leq M |t|,
\end{eqnarray*}
for all $u\in \mathbb{S}^{n-1}$ and $t\in (-\delta_0, \delta_0)$. Then, from the dominated convergence theorem we get
\begin{eqnarray}\label{4.43z}
	\lim_{t\rightarrow 0}\frac{\gamma_n([\tau_t])-\gamma_n(K)}{t}=-\int_{\mathbb{S}^{n-1}}f(v)dS_{0, \gamma_n}(K, v).
\end{eqnarray}
Hence, by \eqref{4.41z} and \eqref{4.43z}, and noting that $K$
is the maximiser of \eqref{4.30z}, we give that for $\lambda\in \mathbb{R}$, 
\begin{eqnarray*}
	0&\geq& \left. \frac{d}{dt}\left(\mathcal{E} \left( [\tau_t]\right)+\lambda\left(\gamma_n([\tau_t])-\frac{1}{2} \right)  \right) \right|_{t=0^+}\\
	&\geq& \int_{\mathbb{S}^{n-1}}f(v)d\mu(v)-\lambda\int_{\mathbb{S}^{n-1}}f(v)dS_{0, \gamma_n}(K, v).
\end{eqnarray*}
Namely, 
\begin{eqnarray}\label{4.44z}
	\int_{\mathbb{S}^{n-1}}f(v)d\mu(v)\leq\lambda\int_{\mathbb{S}^{n-1}}f(v)dS_{0, \gamma_n}(K, v).
\end{eqnarray}
Combining \eqref{4.37z} and \eqref{4.44z}, it follows that 
\begin{eqnarray*}
	\int_{\mathbb{S}^{n-1}}f(v)g(v)dv \leq\frac{\lambda}{\left( \sqrt{2\pi}\right) ^n}\int_{\mathbb{S}^{n-1}} f(v)e^{-\frac{\rho^2_K(\alpha^*_K(v))}{2}}h_K(v)dS(K, v).
\end{eqnarray*}
Since $f$ is arbitrary, we have
\begin{eqnarray}\label{4.45z}
	\frac{dS(K, v)}{dv}\geq \frac{ \left( \sqrt{2\pi}\right) ^ng(v)e^{\frac{\rho^2_K(\alpha^*_K(v))}{2}}}{\lambda h_K(v)},
\end{eqnarray}
and choosing $f\equiv 1$ in \eqref{4.44z} we see
\begin{eqnarray}\label{4.46z}
	\lambda \geq\frac{\mu(\mathbb{S}^{n-1})}{S_{0, \gamma_n}(K, \mathbb{S}^{n-1})}.
\end{eqnarray}
By inequalities \eqref{4.38z} and \eqref{4.45z} we obtain
\begin{eqnarray*}
	\frac{dS(K, v)}{dv}= \frac{ \left( \sqrt{2\pi}\right) ^ng(v)e^{\frac{\rho^2_K(\alpha^*_K(v))}{2}}}{\lambda h_K(v)}
\end{eqnarray*}
for $v\in \mathbb{S}^{n-1}$, where 
\begin{eqnarray*}
	\lambda =\frac{\mu(\mathbb{S}^{n-1})}{S_{0, \gamma_n}(K, \mathbb{S}^{n-1})}
\end{eqnarray*}
from \eqref{4.39z} and \eqref{4.46z}. Hence, $S(K, \cdot)$ is absolutely continuous with respect to the sperical Lebesgue measure. Define
\begin{eqnarray}\label{4.47z}
E_K=\frac{ \left( \sqrt{2\pi}\right) ^nS_{0, \gamma_n}(K, \mathbb{S}^{n-1})g(v)e^{\frac{\rho^2_K}{2}}}{\mu(\mathbb{S}^{n-1}) h_K}.
\end{eqnarray}
Since $K\in \mathcal{K}^n_e$ and $\frac{1}{C}\leq g\leq C$ for some constant $C>0$, we have that there exists a constnat $C_1>0$ such that
\begin{eqnarray}\label{4.48z}
\frac{1}{C_1}\leq E_K\leq C_1.	
\end{eqnarray}

Recalling $\varphi(y)=h_K(y+e)$ where $e\in \mathbb{S}^{n-1}$ and $y\in e^\perp$, it follows from \eqref{3.21z}, \eqref{3.12} and \eqref{3.13} that 
 \begin{eqnarray*}
	\int_\vartheta \mbox{det}D^2 \varphi (y)d\mathcal{H}^{n-1}(y) 
	=\int_{\pi(\vartheta)}\langle v, e\rangle dS(K, v)=\int_{\vartheta}\frac{E_K(\pi(y))}{( \sqrt{1+|y|^2})^{n+1} }d\mathcal{H}^{n-1}(y).
\end{eqnarray*}
Hence, $\varphi$ satisfies the following Monge-Amp\`{e}re equation on $e^\perp$ in the sense of measure:
 \begin{eqnarray*}
 \mbox{det}D^2 \varphi (y)=\frac{E_K(\pi(y))}{( \sqrt{1+|y|^2})^{n+1} }.
\end{eqnarray*}
Together with \eqref{4.48z}, there exists a constant $C_2>0$ such that in a compact set,
 \begin{eqnarray*}
	\frac{1}{C_2}\leq \mbox{det}D^2 \varphi \leq C_2.
\end{eqnarray*}

Arguing in the same way as (i) of Theorem 3.1, it follows from Lemma 3.2 that $\partial K$ is $C^1$ and strictly convex.

{\bf Step 3.}~~Let $d\mu =g(v)dv$ for $v\in \mathbb{S}^{n-1}$ and an even function $g$ with $\frac{1}{C}\leq g\leq C$ where $C>0$. Suppose that $K\in \mathcal{K}^n_e$ is a maximizer to the maximization problem \eqref{4.30z}, which also implies $\gamma_n(K)=\frac{1}{2}$. We next show
\begin{eqnarray*}
	S_{0, \gamma_n}(K, \cdot)=\frac{S_{0, \gamma_n}(K, \mathbb{S}^{n-1})}{\mu(\mathbb{S}^{n-1})}\mu.
\end{eqnarray*}

For any $f\in C(\mathbb{S}^{n-1})$, there is a $\delta>0$ such that $h_K(v)+tf(v)>0$ for all $v \in \mathbb{S}^{n-1}$ and $t\in (-\delta, \delta)$. Then $[h_t]$ is the convex body given in \eqref{4.22z}. Denote by $h(\cdot, t)$ and $\rho (\cdot, t)$ the support function and radial function of $[h_t]$. Define
\begin{eqnarray*}
	\mathscr{L}_2(t, \lambda)=\mathcal{E}([h_t])+\lambda \left( \gamma_n([h_t])-\frac{1}{2}\right).
\end{eqnarray*}
Since $K$ is the maximizer to \eqref{4.30z}, and note also that $\partial K$ is $C^1$ and strictly convex, it follows from \eqref{4.29z}, \eqref{4.23z}, \eqref{4.25z}, \eqref{2.10z} and \eqref{1.1a} that
\begin{eqnarray*}
		0&=& \left.\frac{\partial}{\partial t}\mathscr{L}_2(t, \lambda)\right|_{t=0} \\
	&=&\left.\frac{\partial}{\partial t}\left( \mathcal{E}([h_t])+\lambda \left( \gamma_n([h_t])-\frac{1}{2}\right)\right) \right|_{t=0}\\
&=&\left.\frac{\partial}{\partial t}\left( -\int_{\mathbb{S}^{n-1}}\log h(v, t)d\mu(v)+\lambda \left( \gamma_n([h_t])-\frac{1}{2}\right)\right) \right|_{t=0}\\
&=&-\int_{\mathbb{S}^{n-1}}\frac{f(v)}{h_K(v)}d\mu(v)+\frac{\lambda}{\left(\sqrt{2\pi} \right)^n }\int_{\mathbb{S}^{n-1}}e^{-\frac{\rho_K(u)^2}{2}}\rho_K(u)^n\frac{f(\alpha_K(u))}{h_K(\alpha_K(u))}du
\end{eqnarray*}
\begin{eqnarray*}
&=&-\int_{\mathbb{S}^{n-1}}\frac{f(v)}{h_K(v)}d\mu(v)+\frac{\lambda}{\left(\sqrt{2\pi} \right)^n }\int_{\mathbb{S}^{n-1}}e^{-\frac{|x|^2}{2}}\frac{f(\nu_K(x))}{\langle x, \nu_K(x)\rangle}\langle x, \nu_K(x)\rangle d\mathcal{H}^{n-1}(x)\\
&=&-\int_{\mathbb{S}^{n-1}}\frac{f(v)}{h_K(v)}d\mu(v)+\lambda\int_{\mathbb{S}^{n-1}}\frac{f(v)}{h_K(v)} dS_{0, \gamma_n}(K, v),
\end{eqnarray*}
i.e., 
\begin{eqnarray}\label{4.49z}
\int_{\mathbb{S}^{n-1}}\frac{f(v)}{h_K(v)}d\mu(v)=\lambda\int_{\mathbb{S}^{n-1}}\frac{f(v)}{h_K(v)} dS_{0, \gamma_n}(K, v),
\end{eqnarray}
By the arbitrariness of $f$, we get 
\begin{eqnarray}\label{4.50z}
\mu=\lambda S_{0, \gamma_n}(K, \cdot).
\end{eqnarray}
Let $f=h_K$ in \eqref{4.49z}. Then 
\begin{eqnarray}\label{4.51z}
	\lambda =\frac{\mu(\mathbb{S}^{n-1})}{S_{0, \gamma_n}(K, \mathbb{S}^{n-1})}.
\end{eqnarray}
Substituting \eqref{4.51z} into \eqref{4.50z}, we see
\begin{eqnarray*}
	S_{0, \gamma_n}(K, \cdot)=\frac{S_{0, \gamma_n}(K, \mathbb{S}^{n-1})}{\mu(\mathbb{S}^{n-1})}\mu.
\end{eqnarray*}

{\bf Step 4.}~~At the position we aim to finish the proof for case $p=0$ in Theorem 1.1.

Let $g_j$ be a sequence of positive and smooth functions on $\mathbb{S}^{n-1}$ such that the measures $\mu_j$ with $d\mu_j=g_jd\mathcal{H}^{n-1}\rightarrow \mu$ weakly on $\mathbb{S}^{n-1}$ as $j\rightarrow\infty$, where $\mu$ is the measure given in Theorem 1.1. Note that for each $j$, we see $\frac{1}{C}<g_j<C$ for constant $C>0$. Then it follows from Step 3 that for each $j$ there exists $K_j\in \mathcal{K}^n_e$ such that 
\begin{eqnarray}\label{4.52z}
	S_{0, \gamma_n}(K_j, \cdot)=\frac{S_{0, \gamma_n}(K_j, \mathbb{S}^{n-1})}{\mu_j(\mathbb{S}^{n-1})}\mu_j
\end{eqnarray}
and $\gamma_n(K_j)=\frac{1}{2}.$
Hence, by (iii) of Theorem 3.1 it follows that if $g_j\in C^\alpha(\mathbb{S}^{n-1})$ for $\alpha\in (0, 1)$, then $h_{K_j}\in C^{2, \alpha} (\mathbb{S}^{n-1})$. Therefore, it follows from the standard theory of uniformly elliptic equation that $K_j$ is smooth and uniformly convex, provided $g_j$ is positive and smooth function.

Let $h_{K_j}(v_j)=\max_{v\in\mathbb{S}^{n-1}}h_{K_j}(v)$ for some $v_j\in \mathbb{S}^{n-1}$. Since $\mu$ is not concentrated on any great subspheres of $\mathbb{S}^{n-1}$,  from \eqref{4.18z} we see that there exist $\delta, \zeta>0$ such that for enough large $j$,
\begin{eqnarray*}
	\int_{\{v\in\mathbb{S}^{n-1}:|\langle v, v_j\rangle|>\zeta\}}g_j(v)dv\geq \delta>0.
\end{eqnarray*}
By this, \eqref{1.1}, \eqref{4.19z}, \eqref{4.20z} and \eqref{4.21z} we have
\begin{eqnarray*}
	0&<& \frac{S_{0, \gamma_n}(K_j, \mathbb{S}^{n-1})}{\mu_j(\mathbb{S}^{n-1})}\delta\\
	&\leq& \frac{S_{0, \gamma_n}(K_j, \mathbb{S}^{n-1})}{\mu_j(\mathbb{S}^{n-1})}\int_{\{v\in\mathbb{S}^{n-1}:|\langle v, v_j\rangle|>\zeta\}}g_j(v)dv\\
	&=& \frac{1}{(\sqrt{2\pi})^n}\int_{\{v\in\mathbb{S}^{n-1}:|\langle v, v_j\rangle|>\zeta\}}e^{-\frac{\left|\nabla h_{K_j}(v)\right|^2+h^2_{K_j} (v)}{2}}h_{K_j}(v)\det\left(\nabla^2h_{K_j}(v)+h_{K_j}(v)I\right)dv\\
	&\leq &\frac{1}{(\sqrt{2\pi})^n}\frac{h^{n}_{K_j}(v_j)}{\zeta e^{\frac{\left( \zeta h_{K_j} (v_j)\right)^2}{2}}}\mathcal{H}^{n-1}(\mathbb{S}^{n-1})\rightarrow 0,
\end{eqnarray*}
as $h_{K_j}(v_j)\rightarrow \infty$. We arrive a contradiction. Hence, $K_j$ is uniformly bounded. The uniformly positive lower bound of $h_{K_j}$ is guaranteed by $\gamma_n(K_j)=\frac{1}{2}$ as discussed in Lemma 4.2. Therefore, $h_{K_j}$ is uniformly bounded between two positive constants.

By Blaschke's selection theorem, $K_j$ converges, after passing to a subsequence, to a convex body $K_0\in \mathcal{K}^n_e$. By weak convergence of the $L_0$ Guassian surface area measures and \eqref{4.52z}, it follows that
\begin{eqnarray*}
	S_{0, \gamma_n}(K_0, \cdot)=\frac{S_{0, \gamma_n}(K_0, \mathbb{S}^{n-1})}{\mu(\mathbb{S}^{n-1})}\mu.
\end{eqnarray*}
Moreover, $\gamma_n(K_0)=\frac{1}{2}.$
\hfill $\square$

\section{\bf Existence of asymmetric solutions}
\indent In this section, we will first prove the existence of smooth solution (without symmetry) to the $L_p$ Gaussian Minkowski problem for $p\geq 1$ when the given data is sufficiently smooth and everywhere positive. The key ingredient in the proof is to establish the $L_p$ Gaussian isoperimetric inequality in the following Theorem 5.1. Then an approximation argument will be used to get the existence of asymmetric weak solution to the problem $p\geq 1$.

The following isoperimetric inequality in Gaussian probability space can be obtained by the Ehrhard inequality (see \cite{1z}).

{\it\noindent{\bf Lemma 5.1 }(Gaussian isoperimetric inequality \cite{1z}).~If $K$ is a convex body in $\mathbb{R}^n$, then 
\begin{eqnarray*}
	S_{\gamma_n}(K)\geq \psi\left(\Upsilon^{-1}(\gamma_n(K))\right).
\end{eqnarray*}
Here, $S_{\gamma_n}(K)= S_{\gamma_n}(K, \mathbb{S}^{n-1})$, 
$\psi(t)=e^{-\frac{t^2}{2}}/\sqrt{2\pi}$, and
\begin{eqnarray*}
\Upsilon(x)=\frac{1}{\sqrt{2\pi}}\int_{-\infty}^xe^{-\frac{t^2}{2}}dt.
\end{eqnarray*}}
The following is to establish the $L_p$ Gaussian isoperimetric inequality. It plays a critical role in solving the smooth solution of $L_p$ Gaussian Minkowski problem for $p\geq 1$.

{\it\noindent{\bf Theorem 5.1 }($L_p$ Gaussian isoperimetric inequality).~Let $K\in \mathcal{K}^n_o$. Then for $p\geq 1$,
	\begin{eqnarray}\label{5.1z}
		S_{p, \gamma_n}(K)\geq n\gamma_n(K)\left(\frac{\psi\left( \Upsilon^{-1}(\gamma_n(K))\right) }{n\gamma_n(K)} \right)^p,
	\end{eqnarray}
where $S_{p, \gamma_n}(K)= S_{p, \gamma_n}(K, \mathbb{S}^{n-1})$.}

{\bf Proof.}~~For $K\in \mathcal{K}^n_o$, we know
\begin{eqnarray*}
	S_{p, \gamma_n}(K)=\frac{1}{\left( \sqrt{2\pi}\right) ^n}\int_{\partial K}e^{-\frac{|x|^2}{2}}\langle x, \nu_K(x)\rangle^{1-p}d\mathcal{H}^{n-1}(x).
\end{eqnarray*}
Obviously, 
\begin{eqnarray*}
	S_{1, \gamma_n}(K)=\frac{1}{\left( \sqrt{2\pi}\right) ^n}\int_{\partial K}e^{-\frac{|x|^2}{2}}d\mathcal{H}^{n-1}(x)=S_{\gamma_n}(K).
\end{eqnarray*}
For any $u\in \mathbb{S}^{n-1}$, we have
\begin{eqnarray}\label{5.2z}
\int_{0}^{\rho_K(u)}e^{-\frac{r^2}{2}}r^{n-1}dr\geq e^{-\frac{\rho_K^2(u)}{2}}\int_{0}^{\rho_K(u)}r^{n-1}dr=\frac{1}{n}e^{-\frac{\rho_K^2(u)}{2}}\rho_K^n(u).
\end{eqnarray}
Hence, from \eqref{5.2z} and \eqref{2.10z} it follows that
\begin{eqnarray*}
	\gamma_n(K)&=&\frac{1}{\left(\sqrt{2\pi} \right)^n}\int_{\mathbb{S}^{n-1}}\int_0^{\rho_{K}(u)}e^{-\frac{r^2}{2}}r^{n-1}drdu\\
	&\geq& \frac{1}{n\left(\sqrt{2\pi} \right)^n}\int_{\mathbb{S}^{n-1}}e^{-\frac{\rho_K^2(u)}{2}}\rho_K^n(u)du\\
	&=&\frac{1}{n\left(\sqrt{2\pi} \right)^n}\int_{\partial K}e^{-\frac{|x|^2}{2}}\langle x, \nu_K(x)\rangle d\mathcal{H}^{n-1}(x)=:\widehat{\gamma}_n(K).
\end{eqnarray*}
Namely,
\begin{eqnarray}\label{5.3z}
\widehat{\gamma}_n(K)\leq \gamma_n(K).
\end{eqnarray}
Define 
\begin{eqnarray*}
	d\gamma_n^*(x)=\frac{1}{n\left(\sqrt{2\pi} \right)^n\widehat{\gamma}_n(K)}e^{-\frac{|x|^2}{2}}\langle x, \nu_K(x)\rangle d\mathcal{H}^{n-1}(x),
\end{eqnarray*}
for $x\in \partial K$. Then $\gamma_n^*$ is a probability measure. According to Jensen's inequality, it follows that for $p\geq 1$,
\begin{eqnarray*}
\left(\frac{S_{p, \gamma_n}(K)}{n\widehat{\gamma}_n(K)}\right)^{-\frac{1}{p}}&=&\left( \frac{1}{n\left( \sqrt{2\pi}\right) ^n\widehat{\gamma}_n(K) }\int_{\partial K}e^{-\frac{|x|^2}{2}}\langle x, \nu_K(x)\rangle^{1-p}d\mathcal{H}^{n-1}(x)\right)^{-\frac{1}{p}}\\
&=&\left( \int_{\partial K}\langle x, \nu_K(x)\rangle^{-p}d\gamma_n^*(x)\right)^{-\frac{1}{p}}\\
&\leq& \left( \int_{\partial K}\langle x, \nu_K(x)\rangle^{-1}d\gamma_n^*(x)\right)^{-1}\\
&=& \left(\frac{1}{n\left( \sqrt{2\pi}\right) ^n\widehat{\gamma}_n(K)} \int_{\partial K}e^{-\frac{|x|^2}{2}}d\mathcal{H}^{n-1}(x)\right)^{-1}\\
&=&\left(\frac{S_{\gamma_n}(K)}{n\widehat{\gamma}_n(K)}\right)^{-1}.
\end{eqnarray*}
This together with inequality \eqref{5.3z} has that for $p\geq 1$,
\begin{eqnarray*}
S_{p, \gamma_n}(K)\geq S_{\gamma_n}(K)^p\left(n\widehat{\gamma}_n(K) \right)^{1-p}\geq S_{\gamma_n}(K)^p\left(n\gamma_n(K)\right)^{1-p}.
\end{eqnarray*}
Combining with Lemma 5.1, we see
\begin{eqnarray*}
	S_{p, \gamma_n}(K)\geq \psi\left(\Upsilon^{-1}(\gamma_n(K))\right)^p\left(n\gamma_n(K)\right)^{1-p}.
\end{eqnarray*}
\hfill $\square$

The following inequality is a direct consequence of the $L_p$ Gaussian isoperimetric inequality.

{\it\noindent{\bf Corollary 5.1.}~~Let $K\in \mathcal{K}^n_o$ such that $\gamma_n(K)=\frac{1}{2}$. Then for $p\geq 1$,
	\begin{eqnarray*}
		S_{p, \gamma_n}(K)\geq \left(\frac{1}{\sqrt{2\pi}} \right)^p\left(\frac{n}{2} \right)^{1-p}. 
	\end{eqnarray*}}
The following lemma established in \cite{73} will be needed, which contains a uniform estimate.

 {\it \noindent{\bf Lemma 5.2 \cite{73}.}~Let $\frac{1}{C}<f<C$ for some $C>0$ and $K\in \mathcal{K}_o^n$ with $\gamma_n(K)\geq\frac{1}{2}$.  If $h_K\in C^2(\mathbb{S}^{n-1})$ is a solution of (\ref{1.1}), then there exists a constant $C'>0$ such that
 	\begin{eqnarray*} 
 	\frac{1}{C'}<h_K<C'.
 	\end{eqnarray*}
 }
  
Here, we apply the degree theory to obtain the existence of smooth solution to the $L_p$ Gaussian Minkowski problem for $p\geq 1$. The proof is based on the idea of Huang, Xi and Zhao \cite{71} and Liu \cite{73}. For the degree theory of second-order nonlinear elliptic operators, the reader may wish to consult the work of Li\cite{2z}.

{\it \noindent{\bf Theorem 5.2.}~For $p\geq 1$, suppose that $f$ is a positive smooth function on $\mathbb{S}^{n-1}$ and $|f|_{L_1}<\left(\frac{1}{\sqrt{2\pi}} \right)^p\left(\frac{n}{2} \right)^{1-p}$. Then there exists a smooth convex body $K\in \mathcal{K}_o^n$ with $\gamma_n(K)> \frac{1}{2}$ such that $h_K$ satisfies \eqref{1.1}.
} 

{\bf Proof.}~~Let $f\equiv c_0>0$ small enough and $|c_0|_{L_1}<\left(\frac{1}{\sqrt{2\pi}} \right)^p\left(\frac{n}{2} \right)^{1-p}$. Then it follows from the intermediate value theorem that there exists a constant solution $h_K=s_0$ with $\gamma_n(K)> \frac{1}{2}$. By Lemma 2.3 we see that the solution $K$ is unique. Define $f_t=(1-t)c_0+tf$ where $t\in [-1, 1]$. Then we consider a family of equations on $\mathbb{S}^{n-1}$:
\begin{eqnarray}\label{5.4z}
	\frac{1}{(\sqrt{2\pi})^n}e^{-\frac{|\nabla h|^2+h^2 }{2}}h^{1-p}\det\left(\nabla^2h+hI\right)=f_t.
\end{eqnarray}
Let $0<f\in C^\alpha(\mathbb{S}^{n-1})$. Then there exists a positive constant $C$ such that $\frac{1}{C}<f, c_0< C$. Obviously, $0<f_t\in C^\alpha(\mathbb{S}^{n-1})$,
\begin{eqnarray}\label{5.5z}
	\frac{1}{C}<f_t< C,
\end{eqnarray}
and
\begin{eqnarray}\label{5.6z}
	|f_t|_{L_1}<\left(\frac{1}{\sqrt{2\pi}} \right)^p\left(\frac{n}{2} \right)^{1-p}.
\end{eqnarray}
If $h_t$ solves the equation \eqref{5.4z} and $\gamma_n([h_t])\geq \frac{1}{2}$, then by Lemma 5.2 we see 
\begin{eqnarray}\label{5.7z}
	\frac{1}{C_1}<h_t< C_1,
\end{eqnarray}
where $C_1$ depends on $C$ coming the upper bound of $f$. Moreover, according to the equation \eqref{5.4z} and Theorem 3.1, we can obtain 
\begin{eqnarray}\label{5.8z}
	|h_t|_{C^{2, \alpha}}\leq C_2,
\end{eqnarray}
where the positive constant $C_2$ depends on $|f|_{C^\alpha}$ and the upper bound $C$ of $f$. By the convexity of $h_t$, we also have
\begin{eqnarray}\label{5.9z}
	|\nabla h_t|< C_3,
\end{eqnarray}
where $C_3$ depends $C_1$. Since $h_t\in C^{2, \alpha}$, it follows from \eqref{5.5z}, \eqref{5.7z} and \eqref{5.9z} that for all $t\in [-1, 1]$, the equation \eqref{5.4z} is uniformly elliptic.

Let $C_c^{2, \alpha}(\mathbb{S}^{n-1})$ denote the subset constituting strictly convex functions from $C^{2, \alpha}(\mathbb{S}^{n-1})$ for $\alpha\in (0, 1)$. Now, we construct an open bounded subset $\mathcal{O}$ of $C_c^{2, \alpha}(\mathbb{S}^{n-1})$:
\begin{eqnarray*}
	\mathcal{O}=\left\lbrace h\in C_c^{2, \alpha}(\mathbb{S}^{n-1}): \frac{1}{C_1}<h<C_1, |h|_{C^{2, \alpha}}<C_2, \gamma_n([h])>\frac{1}{2}\right\rbrace. 
\end{eqnarray*}
Define the operators $\mathcal{F}_t: C_c^{2, \alpha}(\mathbb{S}^{n-1})\rightarrow C^{\alpha}(\mathbb{S}^{n-1})$ by
\begin{eqnarray*}
	\mathcal{F}_t(h)=\det\left(\nabla^2h+hI\right)-(\sqrt{2\pi})^nh^{p-1}e^{\frac{|\nabla h|^2+h^2}{2}}f_t.
\end{eqnarray*}
Then the degree $\deg (	\mathcal{F}_t, \mathcal{O}, 0)$ is well-defined if we can prove the zero points of $\mathcal{F}_t$ are not on the boundary of $\mathcal{O}$ for all $t\in[-1, 1]$, i.e., arguing
 \begin{equation*}
	\mathcal{F}_t^{-1}(0)\cap \partial\mathcal{O}=\emptyset,
\end{equation*}
for all $t\in [0, 1]$. We show it by contradiction, and claim that if $h\in \partial \mathcal{O}$, then $\mathcal{F}_t(h)\neq 0$ for all  $t\in[-1, 1]$. Let $\mathcal{F}_t(h)= 0$. Then $h$ solves \eqref{5.4z}. Thus, by \eqref{5.7z} and \eqref{5.8z} this implies that $\gamma_n([h])=\frac{1}{2}$ only. Hence, it follows from Corollary 5.1 that
	\begin{eqnarray*}
	S_{p, \gamma_n}(K)\geq \left(\frac{1}{\sqrt{2\pi}} \right)^p\left(\frac{n}{2} \right)^{1-p}. 
\end{eqnarray*}
This is a contradiction, since if $h$ is a solution to \eqref{5.4z} then 
	\begin{eqnarray*}
	S_{p, \gamma_n}(K)=|f_t|_{L_1}< \left(\frac{1}{\sqrt{2\pi}} \right)^p\left(\frac{n}{2} \right)^{1-p}. 
\end{eqnarray*}

If we can verify $\deg (\mathcal{F}_1, \mathcal{O}, 0)\neq 0$, then there exists at least one solution $h\in \mathcal{O}$ such that $\mathcal{F}_1(h)=0$, which means that $h$ solves the equation \eqref{1.1}, thus achieving our desired result. Note the fact that 
 \begin{equation*}
	\deg(\mathcal{F}_1, \mathcal{O}, 0)=\deg(\mathcal{F}_0, \mathcal{O}, 0),
\end{equation*}
from Proposition 2.2 in Li \cite{2z}. Therefore, we only prove
\begin{equation*}
	\deg(\mathcal{F}_0, \mathcal{O}, 0)\neq 0.
\end{equation*}

Let $\mathcal{L}_{s_0}: C_c^{2, \alpha}(\mathbb{S}^{n-1})\rightarrow C^{\alpha}(\mathbb{S}^{n-1})$
denote the linearized operator of $\mathcal{F}_0$ at $s_0$. It is not hard to compute
\begin{equation*}
	\mathcal{L}_{s_0}(\phi)=s_0^{n-2}\left(\triangle_{\mathbb{S}^{n-1}}\phi+((n-p)-s_0^2)\phi\right).
\end{equation*}
Since spherical Laplacian is a discrete spectrum, we can find a small enough $c_0$ with $|c_0|_{L_1}<\left(\frac{1}{\sqrt{2\pi}} \right)^p\left(\frac{n}{2} \right)^{1-p}$ such that $\mbox{Ker}(\mathcal{L}_{s_0})=\left\lbrace 0 \right\rbrace $. Hence, $\mathcal{L}_{s_0}$ is invertible. By Propositions 2.3 and 2.4 in Li \cite{2z}, and combining with the fact that $s_0$ is the unique solution of $\mathcal{F}_0(h)=0$, we have
\begin{equation*}
	\deg(\mathcal{F}_0, \mathcal{O}, 0)=\deg(\mathcal{L}_{s_0}, \mathcal{O}, 0)\neq 0.
\end{equation*}

Smoothness of solutions to the equation \eqref{1.1} follows directly from the standard theory of uniformly elliptic equations, thus completing the proof of Theorem 5.2.
  \hfill $\square$

We lastly use an appoximation argument to complete the proof of Theorem 1.2.

{\bf Proof of Theorem 1.2.}~Let $f_j$ be a sequence of positive and smooth functions on $\mathbb{S}^{n-1}$ with $|f_j|_{L_1}<\left(\frac{1}{\sqrt{2\pi}} \right)^p\left(\frac{n}{2} \right)^{1-p}$ such that $\{\mu_j\}$, $d\mu_j=f_jd\mathcal{H}^{n-1}$ on $\mathbb{S}^{n-1}$, converges weakly to $\mu$ in Theorem 1.2. It follows from Theorem 5.2 that there is a sequence of smooth convex bodies $K_j$ with $\gamma_n(K_j)>\frac{1}{2}$ such that
\begin{eqnarray}\label{5.10z}
	\frac{1}{(\sqrt{2\pi})^n}h_{K_j}^{1-p}e^{-\frac{|\nabla h_{K_j}|^2+h^2_{K_j} }{2}}\det(\nabla^2h_{K_j}+h_{K_j}I)=f_j.
\end{eqnarray}
This has
\begin{eqnarray}\label{5.11z}
	S_{p, \gamma_n}(K_j, \cdot)=\mu_j.
\end{eqnarray}
Let $h_{K_j}(v_j)=\max_{v\in\mathbb{S}^{n-1}}h_{K_j}(v)$ for some $v_j\in\mathbb{S}^{n-1}$. Since the support of $\mu$ is not contained in a closed hemisphere, it follows by Lemma 5.4 in \cite{73} that there exist positive constant $c$ and $\xi$ such that for large $j$,
\begin{eqnarray}\label{5.12z}
	\int_{\{v\in \mathbb{S}^{n-1}:\langle v, v_j\rangle>\xi\}}f_jdv\geq c>0.
\end{eqnarray}
From (\ref{2.2}), we have that for any $v\in \mathbb{S}^{n-1}$ with $\langle v, v_j\rangle\geq \xi$,
\begin{eqnarray}\label{5.13z}
	h_{K_j}(v)\geq\langle v, v_j\rangle h_{K_j}(v_j)\geq \xi h_{K_j}(v_j).
\end{eqnarray}
Thus, it follows from \eqref{5.10z}, \eqref{5.13z}, \eqref{4.20z} and \eqref{4.21z} that
\begin{eqnarray*}
	0<c&\leq &\int_{\{v\in \mathbb{S}^{n-1}:\langle v, v_j\rangle>\xi\}}f_jdv \\
	&=&\frac{1}{(\sqrt{2\pi})^n}\int_{\{v\in \mathbb{S}^{n-1}:\langle v, v_j\rangle>\xi\}}e^{-\frac{|\nabla h_{K_j}|^2+h^2_{K_j} }{2}}h_{K_j}^{1-p}\det(\nabla^2h_{K_j}+h_{K_j}I)dv\\
	&\leq &	\frac{1}{(\sqrt{2\pi})^n}e^{-\frac{\left(\xi h_{K_j}(v_j)\right)^2 }{2}}h_{K_j}^{n-p}(v_j)\frac{1}{\xi^p}\mathcal{H}^{n-1}(\mathbb{S}^{n-1})\rightarrow 0,
\end{eqnarray*}
as $h_{K_j}(v_j)\rightarrow \infty$. This is a contradiction. Therefore, $K_j$ is uniformly bounded.

We also deduce that $h_{K_j}$ stay uniformly away from zero. If not, there exist $K_j\in \mathcal{K}^n_o$ with $\gamma_n(K_j)>\frac{1}{2}$ and $w_j\in \mathbb{S}^{n-1}$ such that $\min_{w\in\mathbb{S}^{n-1}}h_{K_j}(w)=h_{K_j}(w_j)\rightarrow 0$. According to the definition of support function, we have
\begin{equation*}
	K_j\subset \left\lbrace x\in \mathbb{R}^n: \langle x, w_j\rangle\leq h_{K_j}(w_j)\right\rbrace .
\end{equation*}
Since $K_j$ is uniformly bounded, there exists a centered ball $B_R$ with radius $R>0$ such that
\begin{equation*}
	K_j\subset B_R\cap\left\lbrace x\in \mathbb{R}^n: \langle x, w_j\rangle\leq  h_{K_j}(w_j)\right\rbrace .
\end{equation*}
It is well-known that
\begin{eqnarray*}
	\gamma_n(\mathbb{R}^n)=\frac{1}{\left(\sqrt{2\pi}\right)^n}\int_{\mathbb{R}^n}e^{-\frac{|x|^2}{2}}dx=1.
\end{eqnarray*}
Then 
\begin{eqnarray*}
	\gamma_n(H^-_{w_j})=\frac{1}{\left(\sqrt{2\pi}\right)^n}\int_{H^-_{w_j}\cap B_R}e^{-\frac{|x|^2}{2}}dx+\frac{1}{\left(\sqrt{2\pi}\right)^n}\int_{H^-_{w_j}\setminus B_R}e^{-\frac{|x|^2}{2}}dx=\frac{1}{2},
\end{eqnarray*}
where 
\begin{eqnarray*}
	H^-_{w_j}:=\left\{x\in \mathbb{R}^n: \langle x, w_j\rangle\leq 0 \right\}
\end{eqnarray*}
denotes the closed halfspace with outer normal vector $w_j$. Hence, when $h_{K_j}(w_j)\rightarrow 0$,
\begin{eqnarray*}
	\gamma_n(K_j)\leq\gamma_n\left(B_R\cap\left\lbrace x\in \mathbb{R}^n: \langle x, w_j\rangle\leq h_{K_j}(w_j)\right\rbrace \right)\leq\frac{1}{2},
\end{eqnarray*}
This cannot occur as $\gamma_n(K_j)> \frac{1}{2}.$ Thus, there is a constant $C>0$ such that for all $j$,
\begin{eqnarray*}
	\frac{1}{C}\leq h_{K_j}\leq C.
\end{eqnarray*}

By Blaschke's selection theorem, there exists a convergent subsequence, denoted also by $K_j$, such that $K_j\rightarrow K_0$ in Hausdorff distance. Thus by the weak convergence of $S_{p, \gamma_n}$ and \eqref{5.11z} it follows that 
\begin{eqnarray*}
	S_{p, \gamma_n}(K_0, \cdot)= \mu.
\end{eqnarray*}
Moreover, $	\gamma_n(K_0)>\frac{1}{2}.$
\hfill $\square$

\section{\bf Existence and uniqueness of smooth solutions}
\indent The existence and uniqueness of smooth solutions to the $L_p$ Gaussian Minkowski problem for $p>n$ will be demonstrated in this section, which is inspired by the works of \cite{5z, 6z} and particularly by \cite{62}. The existence part is proved by the continuity method, whereas the uniqueness part is obtained using the strong maximum principle in \cite[Theorem 3.5]{1a}.

We shall require the following estimate.

 {\it \noindent{\bf Lemma 6.1.}~Let $h_K\in C^2(\mathbb{S}^{n-1})$ solves the equation \eqref{1.1} for $p>n$ and a convex body $K$. Then if $\frac{1}{C}<f<C$ for a constant $C>0$, then there exists $C'>0$ such that
	\begin{eqnarray*} 
		\frac{1}{C'}<h_K<C'.
	\end{eqnarray*}
}

{\bf Proof.}~~We first prove the positive lower bound of $h_K$ for a convex body $K$. Assume that $h_K$ attains its minimum at $w_0\in \mathbb{S}^{n-1}$. This means
\begin{eqnarray*} 
\nabla h_K(w_0)=0~~~~~~~~~~\mbox{and}~~~~~~~~~~~~\nabla^2h_K(w_0)\geq 0.
\end{eqnarray*}
Then, according to the equation \eqref{1.1} we have that for $p>n$,
\begin{eqnarray*} 
C>f(w_0)&=&\frac{1}{(\sqrt{2\pi})^n}h^{1-p}_{K}(w_0)e^{-\frac{\left|\nabla h_K(w_0)\right|^2+h^2_K(w_0) }{2}}\det\left(\nabla^2h_K(w_0)+h_K(w_0)I\right)\\
&\geq&\frac{1}{(\sqrt{2\pi})^n}e^{-\frac{h^2_K(w_0) }{2}}h^{n-p}_{K}(w_0)\rightarrow \infty, 
\end{eqnarray*}
as $h_{K}(w_0)\rightarrow 0$. This is a contradiction. Namely, $h_K$ has a positive lower bound. In the same way, at the maximum point of $h_K$ we can obtain that $h_K$ is bounded from above (see also \cite [Lemma 5.1]{73}).\hfill $\square$

When $p> n$, the following lemma shows that the linearized operator of the equation (\ref{1.1}) is invertible at any positive solution. 

{\it \noindent{\bf Lemma 6.2.}~Let $p> n$ and $0<\alpha <1$, and suppose that $f\in C^\alpha (\mathbb{S}^{n-1})$ is a given positive function and $K\in \mathcal{K}_o^n$. If $h=h_K$ solves the equation (\ref{1.1}), then the linearized operator of (\ref{1.1}) at $h$ is invertible.
} 

{\bf Proof.}~For $K\in \mathcal{K}_o^n$, let $h=h_K$ be a solution of (\ref{1.1}), i.e., $h$ solves the equation (\ref{1.1}) for a positive $f\in C^\alpha (\mathbb{S}^{n-1})$. Moreover, we suppose that $h_\varepsilon=he^{\varepsilon\phi}$ for $\phi\in C^{2,\alpha} (\mathbb{S}^{n-1})$ and $\mathcal{L}_h$ denotes the linearized operator of (\ref{1.1}) at $h$. Then we can calculate
\begin{eqnarray}\label{4.1}
	\mathcal{L}_h(\phi)&=&\left. \frac{d}{d\varepsilon}\det\left(\nabla^2h_\varepsilon+h_\varepsilon \delta_{ij}\right) \right|_{\varepsilon=0}-\left.(\sqrt{2\pi})^nf\frac{d}{d\varepsilon}\left(h_\varepsilon^{p-1}e^{\frac{\left|\nabla h_\varepsilon\right|^2+h^2_\varepsilon }{2}}\right)\right|_{\varepsilon=0} \nonumber \\
	&=&	U^{ij}\left(h_{ij}\phi+h_i\phi_j+h_j\phi_i+h\phi_{ij}+h\phi\delta_{ij}\right)\nonumber \\
	&&-(\sqrt{2\pi})^nf\left[(p-1)h^{p-1}e^{\frac{\left|\nabla h\right|^2+h^2 }{2}}\phi+h^{p-1}e^{\frac{\left|\nabla h\right|^2+h^2 }{2}}\left(h_i^2+h^2\right)\phi\right.\nonumber \\
	&&\left.+h^pe^{\frac{\left|\nabla h\right|^2+h^2 }{2}}h_i\phi_i\right],
\end{eqnarray}
where $U^{ij}$ is the cofactor matrix of $\left(h_{ij}+h\delta_{ij}\right)_{(n-1)\times (n-1)}$. Noting the fact that $\frac{1}{C'}<h<C'$ for a constant $C'>0$ by Lemma 6.1 and the assumption that $p> n$, we obtain
\begin{eqnarray}\label{4.2}
	\mathcal{L}_h(1)&=&U^{ij}\left(h_{ij}+h\delta_{ij}\right)-(\sqrt{2\pi})^nfh^{p-1}e^{\frac{\left|\nabla h\right|^2+h^2 }{2}}\left(p-1+h^2_i+h^2\right) \nonumber \\
	&=&(\sqrt{2\pi})^nfh^{p-1}e^{\frac{\left|\nabla h\right|^2+h^2 }{2}}\left(n-p-h^2_i-h^2\right)<0.
\end{eqnarray}
Since $U^{ij}$ is positive definite, at any minimum point of $\phi$ it follows that
\begin{eqnarray}\label{4.3}
V:=hU^{ij}\phi_{ij}+U^{ij}\left(h_i\phi_j+h_j\phi_i\right)-(\sqrt{2\pi})^nfh^{p}e^{\frac{\left|\nabla h\right|^2+h^2 }{2}}h_i\phi_i\geq 0.
\end{eqnarray}
Let $\mathcal{L}_h(\phi)=0$. Then by (\ref{4.1}) and (\ref{4.3}) we have
\begin{eqnarray*}
0=\phi\mathcal{L}_h(1)+V.
\end{eqnarray*}
By the nonnegativity of $V$, this has
\begin{eqnarray*}
	\phi\mathcal{L}_h(1)\leq 0.
\end{eqnarray*}
Thus, it follows from (\ref{4.2}) that $\phi\geq 0$. Similarly, at any maximum point of $\phi$, we have $\phi\leq 0$. Therefore,
\begin{eqnarray*}
	\phi\equiv 0.
\end{eqnarray*}
This implies that $\mathcal{L}_h$ has a trivial kernel. Thus the linearized operator of (\ref{1.1}) at $h$ is invertible. \hfill $\square$

We begin with proving the existence part in Theorem 1.3. Namely, we will prove

{\it \noindent{\bf Theorem 6.1.}~Let $f$ be a positive smooth function on $\mathbb{S}^{n-1}$. Then for $p>n$, there exists a smooth solution $h$ solving the equation \eqref{1.1}.
}
 
{\bf Proof.}~Let $f\equiv c_0>0$. By the intermediate value theorem, there exists a solution $h_K=s_0$ of the equation (\ref{1.1}). Define $f_t=(1-t)c_0+tf$ for $t\in [0, 1]$. We consider the following family of equations
 \begin{eqnarray}\label{4.4}
	\frac{1}{(\sqrt{2\pi})^n}h^{1-p}_{K}e^{-\frac{\left|\nabla h_K\right|^2+h^2_K }{2}}\det\left(\nabla^2h_K+h_KI\right)=f_t.
\end{eqnarray} 
Since $f$ is a positive smooth function on $\mathbb{S}^{n-1}$, there exists a positive constant $C$ such that $\frac{1}{C}<f, c_0< C$. Obviously, $\frac{1}{C}<f_t< C$. Thus by Lemma 6.1 we have that the solutions of equations (\ref{4.4}), denoted by $h_{K_t}$, are bounded from below and above. That is, $\frac{1}{C_1}<h_{K_t}<C_1$ for some positive constant $C_1$. Let $\mathcal{I}$ be the set defined by
 \begin{eqnarray*}
	\mathcal{I}=\left\{ t\in [0, 1]: \mbox{(\ref{4.4}) has a solution}~ h_{K_t} \in C^{2, \alpha}(\mathbb{S}^{n-1}) \right\}.
\end{eqnarray*}  
If we can prove $\mathcal{I}=[0, 1]$, then when $t=1$, there exists a convex body, say   $K_1\in\mathcal{K}_o^n$, with $h_{K_1} \in C^{2, \alpha}(\mathbb{S}^{n-1})$ such that it solves the equation (\ref{1.1}). Hence, we need to show that $\mathcal{I}$ is both open and closed. Clearly, $0\in 	\mathcal{I}$ since $s_0$ solves the equation (\ref{4.4}) with $t=0$. The openness of $\mathcal{I}$ follows from the implicit  function theorem combining with Lemma 6.2, and the closeness of $\mathcal{I}$ follows from Theorem 3.1. 
 
The higher order regularity of the solution to the equation \eqref{1.1} follows from the standard theory of uniformly elliptic equations.\hfill $\square$

We next show the uniqueness part in Theorem 1.3, which is also stated as follows:

{\it \noindent{\bf Theorem 6.2.}~Let $p>n$ and $f$ be a positive smooth function on $\mathbb{S}^{n-1}$. If the functions $h_1$ and $h_2$ are smooth solutions to the equation \eqref{1.1}, then $h_1\equiv h_2$.
}

{\bf Proof.}~~We shall argue by contradiction. Let $h_1, h_2\in C^{\infty}(\mathbb{S}^{n-1})$ be the solutions of equation (\ref{1.1}) for a positive function $f\in C^\infty (\mathbb{S}^{n-1})$. Without loss of generality, we may assume $h_1>h_2$ somewhere on $\mathbb{S}^{n-1}$. Hence, there exists a constant $t> 1$ such that $th_2-h_1\geq 0$ on $\mathbb{S}^{n-1}$ and $th_2-h_1=0$ at some point $v_0\in \mathbb{S}^{n-1}$. Since $p> n$ and $t>1$, we have $t^{n-p}< 1$. Hence for $v\in \mathbb{S}^{n-1}$, we obtain
\begin{eqnarray}\label{5.1}
	\det\left( \nabla^2(th_2)+(th_2)I\right) &=&t^{n-p}(\sqrt{2\pi})^n(th_2)^{p-1}e^{\frac{\left|\nabla h_2\right|^2+h^2_2 }{2}} f \nonumber \\
	&<&t^{n-p}(\sqrt{2\pi})^n(th_2)^{p-1}e^{\frac{\left|\nabla (th_2)\right|^2+(th_2)^2}{2}} f\nonumber \\
	&<&(\sqrt{2\pi})^n(th_2)^{p-1}e^{\frac{\left|\nabla (th_2)\right|^2+(th_2)^2}{2}} f.
\end{eqnarray}
Note that $t>1$ is used in the first inequality. Moreover,
\begin{eqnarray}\label{5.2}
	\det\left( \nabla^2h_1+h_1I\right)=(\sqrt{2\pi})^nh_1^{p-1}e^{\frac{\left|\nabla h_1\right|^2+h_1^2}{2}}f.
\end{eqnarray}
From (\ref{5.1}) and (\ref{5.2}), we have
\begin{eqnarray}\label{5.3}
	&&\det\left( \nabla^2(th_2)+(th_2)I\right)-\det\left( \nabla^2h_1+h_1I\right)\nonumber \\
	&<& (\sqrt{2\pi})^n(th_2)^{p-1}e^{\frac{\left|\nabla (th_2)\right|^2+(th_2)^2}{2}} f-(\sqrt{2\pi})^nh_1^{p-1}e^{\frac{\left|\nabla h_1\right|^2+h_1^2}{2}}f.
\end{eqnarray}
Let $w=th_2-h_1$. Then,
\begin{eqnarray}\label{5.4}
	&&\det\left( \nabla^2(th_2)+(th_2)I\right)-\det\left( \nabla^2h_1+h_1I\right)\nonumber\\
	&=&\int_0^1\frac{d}{d\lambda}\det\left(\nabla^2(h_1+\lambda w)+(h_1+\lambda w)I\right)d\lambda\nonumber\\
	&=&\sum_{i,j=1}^{n-1}\int_0^1\frac{\partial\det\left(\nabla^2(h_1+\lambda w)+(h_1+\lambda w)I\right)}{\partial\left(\nabla_{ij}(h_1+\lambda w)+(h_1+\lambda w)\delta_{ij}\right)}d\lambda\cdot(\nabla_{ij}w+w\delta_{ij})\nonumber\\
	&=&\sum_{i,j=1}^{n-1}\int_0^1U_\lambda^{ij}d\lambda\cdot(\nabla_{ij}w+w\delta_{ij})\nonumber\\
	&=&\sum_{i,j=1}^{n-1}a_{ij}(\nabla_{ij}w+w\delta_{ij}),
\end{eqnarray}
where the coefficient $a_{ij}=\int_0^1U_\lambda^{ij}d\lambda$, and $U_\lambda^{ij}$ is the cofactor matrix of
\begin{eqnarray*}
	M_\lambda:=\nabla^2(h_1+\lambda w)+(h_1+\lambda w)I.
\end{eqnarray*}
In addition,
\begin{eqnarray}\label{5.5}
	&&(\sqrt{2\pi})^n(th_2)^{p-1}e^{\frac{\left|\nabla (th_2)\right|^2+(th_2)^2}{2}} f-(\sqrt{2\pi})^nh_1^{p-1}e^{\frac{\left|\nabla h_1\right|^2+h_1^2}{2}}f \nonumber\\
	&=&\int_0^1\frac{d}{d\lambda}\left((\sqrt{2\pi})^n(h_1+\lambda w)^{p-1}e^{\frac{\left|\nabla (h_1+\lambda w)\right|^2+(h_1+\lambda w)^2}{2}}f \right)d\lambda\nonumber\\
	&=&\sum_{i=1}^{n-1}\int_0^1\frac{\partial}{\partial\left( \nabla_ih_1+\lambda \nabla_i w\right) }\left((\sqrt{2\pi})^n(h_1+\lambda w)^{p-1}e^{\frac{\left|\nabla (h_1+\lambda w)\right|^2+(h_1+\lambda w)^2}{2}}f \right)d\lambda \nonumber\\
	&&\times\nabla_i w+\int_0^1\frac{\partial}{\partial(h_1+\lambda w)}\left((\sqrt{2\pi})^n(h_1+\lambda w)^{p-1}e^{\frac{\left|\nabla (h_1+\lambda w)\right|^2+(h_1+\lambda w)^2}{2}}f \right)d\lambda w \nonumber\\
	&=&\sum_{i=1}^{n-1}b_i\nabla_i w+cw,
\end{eqnarray}
where
\begin{eqnarray*}
	b_i=\int_0^1\frac{\partial}{\partial\left( \nabla_ih_1+\lambda \nabla_i w\right) }\left((\sqrt{2\pi})^n(h_1+\lambda w)^{p-1}e^{\frac{\left|\nabla (h_1+\lambda w)\right|^2+(h_1+\lambda w)^2}{2}}f \right)d\lambda,
\end{eqnarray*}
and
\begin{eqnarray*}
	c=\int_0^1\frac{\partial}{\partial(h_1+\lambda w)}\left((\sqrt{2\pi})^n(h_1+\lambda w)^{p-1}e^{\frac{\left|\nabla (h_1+\lambda w)\right|^2+(h_1+\lambda w)^2}{2}}f \right)d\lambda.
\end{eqnarray*}
Thus, it follows from (\ref{5.3}), (\ref{5.4}) and (\ref{5.5}) that
\begin{eqnarray*}
	Lw:=\sum_{i,j=1}^{n-1}a_{ij}\nabla_{ij}w-\sum_{i=1}^{n-1}b_iw_i+\left(\sum_{i=1}^{n-1}a_{ii}-c\right)w< 0.
\end{eqnarray*}
Since $h_1, h_2$ are the smooth solutions of equation (\ref{1.1}) with $0<f\in C^{\infty}(\mathbb{S}^{n-1})$, the matrices $M_\lambda$ are smooth, positive definitive for all $\lambda\in [0, 1]$. Thus we have
\begin{eqnarray*}
	\frac{1}{C}I\leq \{a_{ij}\}\leq CI
\end{eqnarray*}
for some positive constant $C$. Namely, $L$ is uniformly elliptic. From the strong maximum principle (see \cite[Theorem 3.5]{1a}), we get $w\equiv0$ on $\mathbb{S}^{n-1}$. In other words, $th_2\equiv h_1$ on $\mathbb{S}^{n-1}$. This implies that $th_1$ solves the equation (\ref{1.1}), i.e.,
\begin{eqnarray*}
	\det\left( \nabla^2(th_2)+(th_2)I\right) =(\sqrt{2\pi})^n(th_2)^{p-1}e^{\frac{\left|\nabla (th_2)\right|^2+(th_2)^2}{2}} f.
\end{eqnarray*}
By this, we see
\begin{eqnarray*}
	\det\left( \nabla^2h_2+h_2I\right) &=&t^{p-n}(\sqrt{2\pi})^nh_2^{p-1}e^{\frac{\left|\nabla (th_2)\right|^2+(th_2)^2}{2}} f \\
	&>& t^{p-n}(\sqrt{2\pi})^nh_2^{p-1}e^{\frac{\left|\nabla h_2\right|^2+h_2^2}{2}} f
	\\	&=&t^{p-n}\det\left( \nabla^2h_2+h_2I\right),
\end{eqnarray*}
which deduce $t< 1$ since $p> n$. This is a contradiction to $t>1$.
\hfill $\square$

\vskip10pt

\bibliographystyle{amsplain}

\end{document}